\documentclass[12pt]{article}
\usepackage{graphicx} 
\usepackage{amsmath} 
\usepackage{amssymb} 
\usepackage{amsfonts} 
\usepackage{amsthm} 
\usepackage{longtable} 
\usepackage{comment}

\usepackage{color}
\usepackage{subfig}
\input{epsf}
\input xy
\xyoption{all}
\usepackage{latexsym}
\newtheorem{lemma}{Lemma}[section] 
\newtheorem{theorem}[lemma]{Theorem} 
\newtheorem{proposition}[lemma]{Proposition}

\newtheorem{example}[lemma]{Example} 
\newtheorem{definition}[lemma]{Definition} 
   
\renewcommand{\epsilon}{\varepsilon} 
 
\renewcommand{\hat}{\widehat} 

\renewcommand{\leq}{\leqslant} 
\renewcommand{\geq}{\geqslant}

\newcommand{\Po}{\mathcal{P}}
\newcommand{\R}{\mathcal{R}}
\newcommand{\Q}{\mathcal{Q}} 
\newcommand{\F}{\mathcal{F}} 
\newcommand{\T}{\mathcal{T}}

\newcommand{\N}{\mathrm{N}}

\newcommand{\K}{\mathcal{K}}

\newcommand{\Ga}{\Gamma}

\newcommand{\tom}{\mathcal{T}_{hh}}
\newcommand{\sem}{semiregular\,\,}

\newcommand{\rk}{\mathrm{rank}} 
\begin{document}


\title{Semiregular Polytopes and Amalgamated C-groups}
 
\author{B. Monson\thanks{Supported by NSERC of Canada Discovery Grant \# 4818}\\
University of New Brunswick\\
Fredericton, New Brunswick, Canada E3B 5A3
\and and\\[.05in]
Egon Schulte\thanks{Supported by NSF-grant DMS--0856675}\\
Northeastern University\\
Boston, Massachussetts,  USA, 02115}

\date{ \today }
\maketitle

\begin{abstract}
In the classical setting,  a convex polytope is said to  be semiregular if its facets are regular and its symmetry group is transitive on vertices. This paper studies semiregular abstract polytopes, which have abstract regular facets, still with combinatorial automorphism group transitive on vertices. We analyze the structure of the automorphism group, focusing in particular on polytopes with two kinds of regular facets occurring in an ``alternating" fashion. In particular we use group amalgamations to prove that given two compatible $n$-polytopes $\Po$ and $\Q$, there exists a universal abstract semiregular $(n+1)$-polytope which is obtained by ``freely" assembling alternate copies of $\Po$ and $\Q$. We also employ  modular reduction techniques to construct finite semiregular 
polytopes from reflection groups over finite fields.  



\end{abstract}

\section{Introduction}
In the classical setting,  a convex $n$-polytope $\Po$ is said to be 
\textit{uniform} if its facets are uniform and its symmetry group 
is transitive on vertices  (see Coxeter~\cite{wyt1,crsrp1}, Johnson~\cite{johns}). To start this inductive definition
in a pleasant way, we agree that uniform polygons should be regular.
The same definition can be transferred to the abstract  
(combinatorial) setting. But  \textit{all} polygons are combinatorially  regular,
so  one  soon suspects that the  abstract uniform polytopes  form a 
huge, perhaps untamable class of strange objects. 
Indeed, the abstract  uniform polytopes  
generalize the abstract regular polytopes, which
in a sense are `maximally' symmetric and which are already quite abundant. 

In this paper,
we will  generalize regularity more modestly and
focus on  \textit{semiregular} polytopes $\Po$,
which have regular facets, still with automorphism group 
$\Gamma(\Po)$ transitive on vertices (see  Coxeter~\cite{crsrp1,crsrp2,crsrp3}).
All of the classical Archimedean solids, for example,
are convex semiregular $3$-polytopes. 
Perhaps the best-known  non-regular example is the 
\textit{cuboctahedron}, which  can be
obtained either by truncating a cube to its edge midpoints,
or by assembling squares and equilateral triangles, two of each placed alternately 
around each vertex. 

This  sort of behaviour also appears in the
familiar  tiling $\T$ of 
Euclidean $3$-space  by regular octahedra 
and tetrahedra, the beginning  of which is displayed in Figure~\ref{fig1}. 
In fact, $\T$ is an infinite semiregular $4$-polytope.
Our main concern in this paper will be   
abstract semiregular polytopes like this, with two kinds of regular facets
occurring in  an `alternating' fashion. The essential features of
our construction are contained in Theorem~\ref{wytsemi}.

Later, in Theorem~\ref{univ}, we  prove the existence of a universal
abstract semiregular (n+1)-polytope $\mathcal{U}_{\Po, \Q}$,
which is obtained by ``freely" assembling alternate copies
of two compatible $n$-polytopes $\Po$ and $\Q$. Finally, in Section~\ref{finfi}, 
we employ  modular reduction techniques to construct finite semiregular 
polytopes from reflection groups over finite fields.  

\section{Abstract Polytopes and their Automorphism Groups}\label{abst}

An abstract  $n$-polytope $\Po$  has some of the key 
combinatorial properties of the face lattice of a convex $n$-polytope;
in general, however, $\Po$ need not be a lattice, 
need not be finite, need not have any familiar geometric realization.
Let us  summarize some general definitions and results, referring  to McMullen \& Schulte~\cite{arp}
for details.
An \textit{abstract $n$-polytope} $\Po$ is a partially ordered
set with  properties \textbf{A}, \textbf{B} and 
\textbf{C} below.
\medskip

\noindent\textbf{A}: $\Po$ has a strictly monotone rank function with range 
$\{-1,0,\ldots,n\}$. 

\medskip\noindent
An element $F \in \Po$ with $\rk(F)=j$ is
called a $j$-\textit{face}; often  $F_j$ will indicate a $j$-face. Moreover,
$\Po$ has a unique least face
$F_{-1}$ and   unique greatest face $F_n$. 
Each maximal chain or \textit{flag} in $\Po$ therefore contains $n+2$ faces,
so that $n$ is the number of \textit{proper} faces in each flag.
We let $\F(\Po)$ be the set of all flags in $\Po$.
Naturally, faces of ranks 0, 1 and $n-1$ are called 
vertices, edges and  facets, respectively.
\medskip 

\noindent\textbf{B}: Whenever $F < G$ with $\rk(F)=j-1$ and
$\rk(G)=j+1$, there are exactly two $j$-faces $H$ with
$F<H<G$.

\medskip\noindent
For $0 \leq j \leq n-1$ and any flag $\Phi$, 
there thus exists a unique \textit{adjacent} 
flag $\Phi^j$, differing from $\Phi$ in just the
face of rank $j$ . With this notion of adjacency,  $\F(\Po)$
becomes the \textit{flag graph} for $\Po$.
If $F \leq G$ are incident faces  
in  $\Po$,  we call
$$ G/F := \{ H \in \mathcal{P}\, | \, F \leq H \leq G \}\;.$$
a \textit{section}  of $\Po$.
 
\medskip
\noindent\textbf{C}: $\Po$ is \textit{strongly flag--connected}, 
that is, the flag graph for each section is connected. 
 
\medskip\noindent
It follows that 
$G/F$ is a ($k-j-1$)-polytope in its own right, if 
$F \leq G$ with $\mathrm{rank}(F) = j \leq k = \mathrm{rank}(G)$.
In particular, $F_n/F$ is the {\em co-face\/} of a face $F$ in $\Po$, and its
rank $n-1-\mathrm{rank}(F)$ is the {\em co-rank} of $F$.
For example, if $F$ is a vertex, then the
section $F_n/F$ is called the \textit{vertex-figure} over $F$.
Likewise, it is useful to think of the  $k$-face $G$ 
as having the structure of the $k$-polytope $G/F_{-1}$.

The \textit{automorphism group} $\Gamma(\Po)$ consists of all
order-preserving bijections on $\Po$. We say  $\Po$ is \textit{regular} if 
$\Gamma(\Po)$ is transitive on the flag set $\F(\Po)$. In this case we
may choose any one flag $\Phi \in \F(\Po)$ as \textit{base flag},  
then  define $\rho_j $ to be  the (unique) automorphism  
mapping $\Phi$ to $\Phi^j$, for $0 \leq j \leq n-1$. Each $\rho_j$ has period $2$.
 From \cite[2B]{arp}
we recall that $\Gamma(\Po)$ is then a \textit{string C-group}, meaning that 
it  has  the following properties \textbf{SC1} and \textbf{SC2}:

\medskip\noindent\textbf{SC1}:   $\Gamma(\Po)$ is 
generated by $ \{\rho_0, \ldots, \rho_{n-1}\}$. These  involutory generators
satisfy the commutativity relations
typical of a Coxeter group with string diagram, namely  
\begin{equation}\label{relreg}
(\rho_ j \rho_k)^{p_{jk}} = 1,\;  \mathrm{ for }\; 
0 \leq j \leq k \leq n-1, 
\end{equation}
where   $p_{jj} = 1$ and $p_{jk} = 2$ whenever $|j-k|>1$.

\medskip\noindent\textbf{SC2}: $\Gamma(\Po)$ satisfies the
\textit{intersection condition}
\begin{equation}\label{interreg}
\langle I\rangle  \cap \langle J\rangle  = \langle I \cap J\rangle, 
\;  \mathrm{ for\,\, any }\; I, J \subseteq \{\rho_0, \ldots, \rho_{n-1}\}\;.
\end{equation}

\medskip

\noindent
The fact that one can reconstruct 
a regular polytope in a canonical way 
from any string C-group $\Gamma$ is at the heart of the theory \cite[2E]{arp}.

The periods $p_j := p_{j-1,j}$ in (\ref{relreg}) 
satisfy $2 \leq p_j \leq \infty$ and        
are assembled into
the \textit{Schl\"{a}fli} symbol $\{ p_1, \ldots, p_{n-1}\}$ for the
regular polytope $\Po$. We note again that every $2$-polytope or \textit{polygon} 
$\{p_1\}$ is automatically  abstractly regular; its automorphism group
is the dihedral group $\mathbb{D}_{2 p_1}$ of order $2 p_1$. 

There are various ways to relax symmetry and thereby
broaden the class  of groups $\Gamma(\Po)$. As we suggested in \S1, the 
class of uniform polytopes is far too broad for our purposes here.
Thus we restrict our reach somewhat with

\begin{definition} An abstract polytope $\Po$ is \emph{semiregular}
if it has regular facets and its automorphism group $\Gamma(\Po)$ is transitive on vertices
\emph{\cite[p. 77]{martini}}.
\end{definition}

Every regular polytope is clearly semiregular.
In this paper  we focus our investigation of \sem polytopes on a particularly interesting
subclass derived from Wythoff's construction described in the next section. We might call
these \textit{alternating semiregular polytopes}, since
they have facets of possibly two distinct types appearing in
alternating fashion around faces of co-rank $2$.

\section{A Concrete (Geometrical) Version of Wythoff's Construction}
\label{gcwyt}
\begin{figure}[htbp]
\centering
\includegraphics[width=90mm]{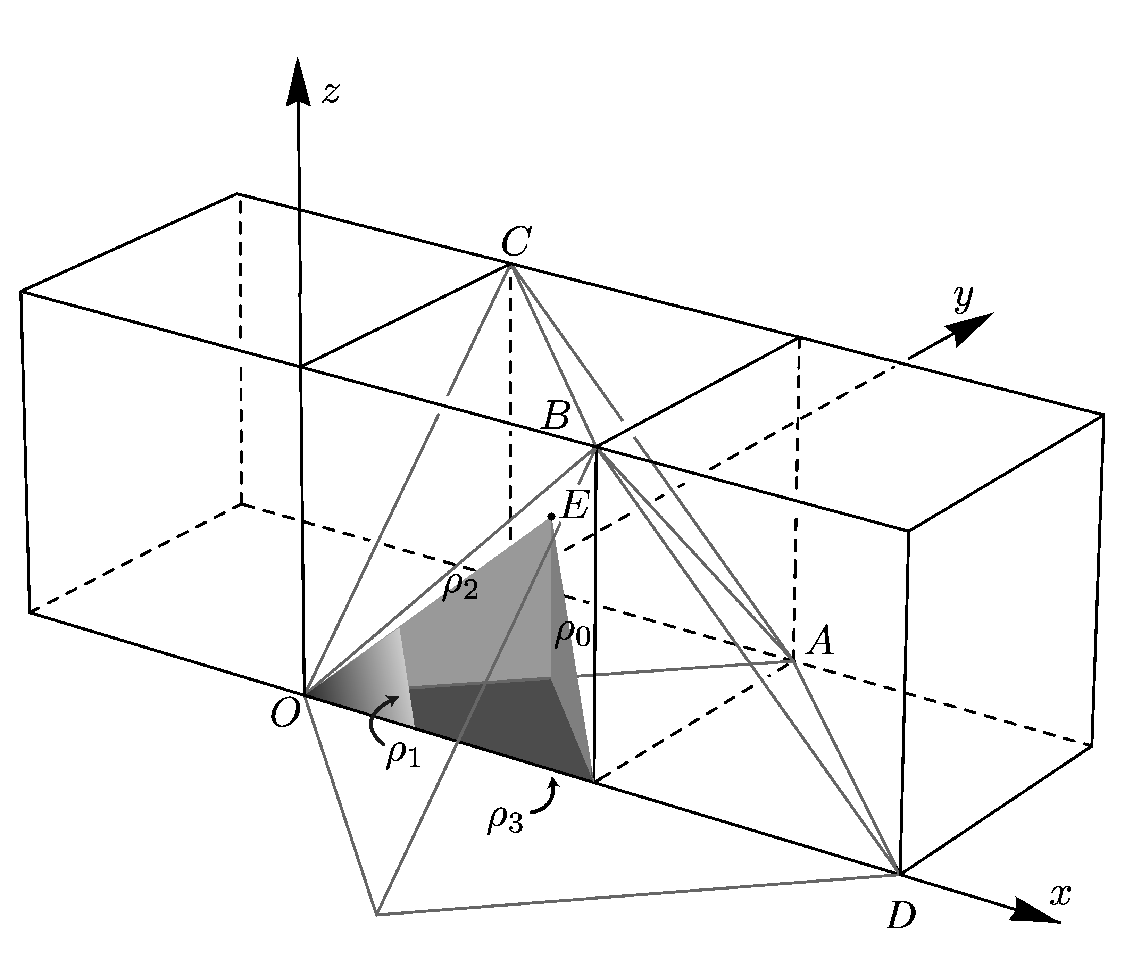}
\caption{The birth of the $4$-polytope $\T$, a \sem
tessellation of $\mathbb{R}^3$.}\label{fig1}
\end{figure}

Let us take a closer look at the uniform tessellation $\T$ mentioned earlier.
As a combinatorial object, this tessellation of $\mathbb{R}^3$
is an abstract \sem $4$-polytope whose combinatorial automorphism
group  can be identified  
 with its geometric symmetry group (see Figure~\ref{fig1}).

A simple way to describe $\T$ is to  first   
imagine  Euclidean space $\mathbb{R}^3$  tiled as usual by unit cubes.
(Although this tiling is itself a regular $4$-polytope, for now it will serve
mainly as scaffolding for $\T$.)   Each cube has two inscribed 
regular tetrahedra. Pick one in each cube, starting 
with the tetrahedron with vertices
$O = (0,0,0), A = (1,1,0), B = (1,0,1), C = (0,1,1)$  
in the standard unit cube, then alternating thereafter as one passes between
adjacent cubes. 
We thus get the tetrahedral facets of $\T$; the octahedral
facets tile what is left of  $\mathbb{R}^3$. Every vertex of $\T$ is 
surrounded by six octahedra and eight tetrahedra; indeed,
each vertex-figure
is a cuboctahedron.
(Compare \cite[\S 4.7]{RP}, where Coxeter describes $\T$ as a 
\textit{quasiregular} tessellation, with modified
Schl\"{a}fli symbol  $\{3, \hspace*{1mm}_4^3 \}$. The usage of the term ``quasiregular" in \cite{RP} 
implies the local alternating behaviour we focus on in this paper.)

Notice that the Euclidean symmetry group $\Gamma(\T)$ contains the
face-centred cubic lattice generated by translations
$\tau_1, \tau_2, \tau_3$ along the edges $OA, OB, OC$
of the base tetrahedron.  The point group stabilizing vertex $O$
is the octahedral group of order $48$ generated by  reflections
$\rho_1, \rho_2, \rho_3$, whose mirrors 
are indicated in
Figure~\ref{fig1}, along with the mirror
for a fourth reflection $\rho_0$. These four
reflections generate $\Gamma(\T)$, and their mirrors
enclose a tetrahedral fundamental region for the action
of $\Gamma(\T)$ on $\mathbb{R}^3$. In fact,
$\Gamma(\T)$
is the (infinite) Coxeter group of type $\widetilde{B}_3$. 

It is convenient to blur the distinction between the affine reflections
$\rho_j$ and their abstract counterparts in a presentation of   
$\Gamma(\T)$ as the affine  Coxeter group $\widetilde{B}_3$ with Coxeter diagram  
\begin{equation}
\label{fig2}
\centering
\begin{picture}(180,25)
\put(25,0){
\multiput(15,0)(45,0){2}{\circle*{5}}
\multiput(95,-23)(0,46){2}{\circle*{5}}
\put(15,0){\circle{9}}
\put(15,0){\line(1,0){45}}
\put(60,0){\line(3,2){35}}
\put(60,0){\line(3,-2){35}}
\put(74,15){\scriptsize{$4$}}
\put(12,8){\scriptsize $\rho_0$}
\put(56,8){\scriptsize $\rho_1$}
\put(100,-25.5){\scriptsize $\rho_2$}
\put(100,21.5){\scriptsize $\rho_3$}}
\end{picture}
\end{equation}
\vskip.37in
\noindent
Recall from \cite[\S6.5]{humph} that $\Gamma(\T)$ has these defining relations
on its standard generators:
\begin{equation}\label{B3cox}
\rho_j^2 = (\rho_0 \rho_2)^2 = (\rho_0 \rho_3)^2 = (\rho_2 \rho_3)^2
= (\rho_0 \rho_1)^3 = (\rho_1 \rho_2)^3 = (\rho_1 \rho_3)^4 = 1\; .
\end{equation}
Evidently $\Gamma(\T)$ is not a string C-group. Indeed, $\T$ is not regular,
for it has two kinds of facets and two flag orbits. Nevertheless,
$\Gamma(\T)$   does satisfy the intersection
condition (\ref{interreg}),
as does any Coxeter group \cite[Theorem 5.5]{humph}. This fact
reappears in Definition~\ref{defS} below.



The ringed node of the diagram in (\ref{fig2}) is an ingenious decoration  invented by 
Coxeter \cite{wyt1} and is meant to  encode 
\textit{Wythoff's construction} for $\T$. The essential idea is that each 
$j$-face of $\T$ lies in the same $\Gamma(\T)$-orbit as a special
$j$-face $F_j$, whose stabilizer $\Sigma(F_j)$ is a certain parabolic
subgroup of $\Gamma(\T)$.  This $\Sigma(F_j)$  is generated by the $\rho_k$'s 
corresponding to nodes in a subdiagram, which in turn consists of
a connected `active' part, which has $j$ nodes 
including the
ringed node, and   a `passive' part induced on all nodes not connected to
the active part. For example, there is one base vertex
$F_0 = O$ in Figure~\ref{fig1}; it has empty active part 
and is fixed (passively) by 
$\Sigma(F_0) = \langle \rho_1, \rho_2, \rho_3 \rangle$. 
The base edge
$F_1$ has vertices $O$ and $A = (O)\rho_0$, the image  of $O$ under $\rho_0$; the active part of
$\Sigma(F_1) = \langle \rho_0, \rho_2, \rho_3 \rangle$
is $\langle \rho_0 \rangle$. The base equilateral triangle
$F_2$ has vertices $OAB$ and $\Sigma(F_2) = \langle \rho_0, \rho_1 \rangle$
is totally active.  Finally, there are  two special $3$-faces, $F_{3}=F_{3}^{(1)}, F_{3}^{(2)}$ (say),
namely an octahedral facet
$F_3^{(1)}$ with $\Sigma(F_3^{(1)}) = \langle \rho_0, \rho_1, \rho_3 \rangle$,
and a tetrahedral facet $F_3^{(2)}$ with 
$\Sigma(F_3^{(2)}) = \langle \rho_0, \rho_1, \rho_2 \rangle$.
In this `geometrical' version of Wythoff's construction, 
each basic face $F_j$ is the convex
hull of the $\Sigma(F_j)$-orbit of $O$.
Observe that  general $j$-faces in the  $\Gamma(\T)$-orbit of
$F_j$ correspond exactly to  the right cosets of
 $\Sigma(F_j)$ in $\Gamma(\T)$.

We note that the description above can be adjusted somewhat 
to accommodate other sorts of Coxeter diagrams, with arbitrary 
sets of ringed nodes \cite{wyt1}. 
(Active parts may then become disconnected.)     
However, the regular case works as expected:     
just ring one terminal node in a string diagram \cite[1B]{arp}. 
Consider, for example, this subdiagram extracted 
from the diagram in (\ref{fig2}): 

\begin{picture}(180,72)
\put(147,35){
\multiput(15,0)(45,0){2}{\circle*{5}}
\put(95,-23){\circle*{5}}
\put(15,0){\circle{9}}
\put(15,0){\line(1,0){45}}
\put(60,0){\line(3,-2){35}}
\put(12,8){\scriptsize $\rho_0$}
\put(56,8){\scriptsize $\rho_1$}
\put(100,-25.5){\scriptsize $\rho_2$}}
\end{picture}

\noindent
The subgroup $\langle \rho_0, \rho_1, \rho_2\rangle$ of $\widetilde{B}_3$
is, in fact, the $A_3$ Coxeter group (isomorphic to the symmetric  group
$\mathbb{S}_4$); and Wythoff's construction now produces
a  regular tetrahedron, such as  $OABC$ in Figure~\ref{fig1}.

Let us summarize the discussion above in

\begin{example}\label{tessT}\emph{The $4$-polytope $\T$ defined
by the diagram in (\ref{fig2}) is a  semiregular
tessellation of $\mathbb{R}^3$ by regular octahedra and tetrahedra.
Its symmetry group $\Gamma(\T)$ is the Coxeter group $\widetilde{B}_3$. 
}\end{example}

We refer to \cite{chen, modpat} for a 
much broader look at the concrete geometrical aspects
of Wythoff's construction. Below we shall pursue instead
an abstract, that is to say, combinatorial  generalization of the   
construction,
motivated   by the above example, though   still tailored to our immediate
needs.

\section{An Abstract (Combinatorial) Version of Wythoff's Construction}
\label{combwyt}
Suppose 
that $\Ga = \langle \alpha_0, \ldots, \alpha_{n-2},\alpha_{n-1},\beta_{n-1} 
\rangle$ is a group generated by involutions which satisfy 
the commutativity relations
implicit in the {\em tail-triangle diagram\/} 
\begin{equation}
\label{2kgp}
\centering
\begin{picture}(180,30)
\put(100,0){
\multiput(15,0)(45,0){2}{\circle*{5}}
\multiput(100,-27)(0,53.3){2}{\circle*{5}}
\multiput(-130,0)(45,0){2}{\circle*{5}}
\put(-130,0){\circle{9}}
\put(-5,0){\line(1,0){20}}
\put(-85,0){\line(1,0){20}}
\put(-130,0){\line(1,0){45}}
\put(15,0){\line(1,0){45}}
\put(60,0){\line(3,2){40}}
\put(60,0){\line(3,-2){40}}
\put(100,-27){\line(0,1){53.3}}
\put(4,9){\scriptsize $\alpha_{n-3}$}
\put(48,9){\scriptsize $\alpha_{n-2}$}
\put(-134,9){\scriptsize $\alpha_{0}$}
\put(-89,9){\scriptsize $\alpha_{1}$}
\put(-49,-0.5){$\ldots\ldots$}
\put(105,-29){\scriptsize $\beta_{n-1}$}
\put(105,25){\scriptsize $\alpha_{n-1}$}
\put(103.5,-3){\scriptsize k}}
\end{picture}
\end{equation}
\vskip.3in
\noindent
The label `$k$' indicates that $\alpha_{n-1} \beta_{n-1}$ has 
period $k$, for some $k=2, \ldots, \infty$.  However,
all other periods of products of two `adjacent' generators are unspecified
for the moment. Then the group $\Ga$ is called a {\em tail-triangle group}. 
We have already encountered such a diagram in 
(\ref{fig2}), in which $k=2$ for the non-adjacent nodes labelled 
$\rho_2$ and $\rho_3$. When $n=2$ there is no ``tail" and $\Ga$ is 
{\em triangle group\/} with a {\em triangle diagram}; this is consistent with 
standard terminology (see \cite{jones1}). We allow the degenerate (base) case 
$n=1$ when $\Ga = \langle\alpha_0,\beta_0\rangle$ is just the 
dihedral group $\mathbb{D}_{2k}$.   

Suppose also
that $\Ga$ satisfies the 
intersection condition on its distinguished subgroups.
In other words, for all subsets 
$I, J  \subseteq \{\alpha_0, \ldots, \alpha_{n-2},\alpha_{n-1},\beta_{n-1}\}$,
we assume that
\begin{equation}\label{intercon}
\langle I\rangle \cap \langle J\rangle  = \langle I\cap J \rangle\,.
\end{equation}
Thus $\Ga$ is a  C-group (see \cite[2E]{arp}). 

It follows that the subgroups 
$$\Ga_n^{\Po} : = \langle \alpha_0, \ldots, \alpha_{n-2},\alpha_{n-1}\rangle$$
and
$$\Ga_n^{\Q} : = \langle \alpha_0, \ldots, \alpha_{n-2},\beta_{n-1} \rangle$$
are string C-groups, indeed automorphism groups for 
regular $n$-polytopes $\Po$ and $\Q$, respectively.

\begin{definition} 
We call a group $\Ga = \langle \alpha_0, \ldots, \alpha_{n-2},\alpha_{n-1},\beta_{n-1} 
\rangle$ represented by a tail-triangle diagram as in 
(\ref{2kgp}) and satisfying the intersection property (\ref{intercon}) a {\em tail-triangle C-group\/}. 
\end{definition}
We can now describe a 
combinatorial version of the Wythoff construction implied by ringing the node
labelled $\alpha_0$ in (\ref{2kgp}). The actual details follow fairly closely
those for the regular case described  in \cite[2E]{arp}.
Anticipating Theorem~\ref{wytsemi}, 
let us denote the resulting $(n+1)$-polytope $\mathcal{S}$,
or $\mathcal{S}(\Gamma)$, if we wish to emphasize the underlying group. Keep in mind
that any such group comes with a specified list of generators, arranged as 
in (\ref{2kgp}).

\begin{definition}\label{defS}
\emph{\textbf{The \sem (n+1)-polytope $\mathcal{S} = \mathcal{S}(\Ga)$}.}

Suppose that the group
$\Ga = \langle \alpha_0, \ldots, \alpha_{n-2},\alpha_{n-1},\beta_{n-1} \rangle$
is a tail-triangle C-group.
Take the \emph{improper} faces of $\mathcal{S}$ to be
two distinct copies $\Ga_{-1}$ and $\Ga_{n+1}$ of $\Ga$. Next, for 
$0 \leq j \leq n-2$,  define the \emph{$j$-faces} of $\mathcal{S}$ to be 
all right cosets in $\Ga$ of 
$$\Ga_j := \langle \alpha_0, \ldots,\alpha_{j-1}, \alpha_{j+1},\ldots,
  \alpha_{n-1},\beta_{n-1}\rangle\;.$$
The $(n-1)$-faces, or \emph{ridges}, of $\mathcal{S}$ are all right cosets of 
$$\Ga_{n-1} := \langle \alpha_0, \ldots, \alpha_{n-2}\rangle\;.$$
Finally, the $n$-faces, or \emph{facets}, of $\mathcal{S}$ are all right cosets of 
either $\Ga_n^{\Po}$ or $\Ga_n^{\Q}$. When there is no risk of confusion
we use $\Ga_n$ to denote \emph{either} of these subgroups.

Finally we define what turns out  to be a partial order on
the set of all such faces by taking
$$
 \Ga_j \nu <  \Ga_k \mu
$$
whenever $-1 \leq j < k \leq n+1$ and 
$ \Ga_j\nu \cap \Ga_k\mu \neq \emptyset$, for $\mu,\nu \in \Ga$. 
\end{definition}

\noindent\textbf{Remarks}.
There is a  reason for indicating strict inequality `$<$' here.
For $j = k = n$, we must explicitly  forbid cosets
$\Ga_n^{\Po}\mu$ and  $\Ga_n^{\Q}\nu$ from being incident.
(Occasionally these are non-disjoint; but distinct facets must never
 be incident. But for that   
glitch, we could say that $\mathcal{S}$ is 
  a thin coset geometry.) 
Anyway,  considering this we naturally agree that 
$ \Ga_j \nu \leq  \Ga_k \mu$  if and only if
$ \Ga_j \nu <  \Ga_k \mu$ or (when $j = k$) $ \Ga_j \nu = \Ga_k \mu$
(cf. Lemma~\ref{dist}(b) below). 

Before we move on, note that (in contrast to 
standard notation for string C-groups) 
$\Ga_j$ is here not always equal to
the subgroup of $\Ga$ generated by all generators 
but the $j$th (if applicable).    

We will now 
prove that $\mathcal{S}$ is an abstract  polytope. We begin with a few 
lemmas  concerning 
a tail-triangle group $\Ga$.

\begin{lemma}\label{dist}
{\rm (a)} The  distinguished subgroups $\langle J \rangle$ are pairwise 
distinct. In particular, the subgroups $\Ga_0, \ldots, \Ga_{n-1},
\Ga_n^{\Po}, \Ga_n^{\Q}$ are mutually distinct and never equal $\Ga$ itself.

{\rm (b)} Suppose   $\Ga_j \mu = \Ga_k\nu$, for $\mu ,\nu \in \Ga$
and $0 \leq j \leq k \leq n$. Then 
$j = k$  and $\Ga_j = \Ga_k$. Furthermore, 
cosets $\Ga_n^{\Po} \mu$ and $ \Ga_n^{\Q} \nu$ can never be equal.

{\rm (c)} Let $0 \leq j_1 < j_2 < \ldots  <j_m \leq n$.
Suppose that $\mu_{j_1}, \ldots, \mu_{j_m} \in \Ga$ with 
$$\Ga_{j_i} \mu_{j_i} \cap \Ga_{j_{i+1}} \mu_{j_{i+1}}\neq \emptyset \; ,$$ 
for $1\leq i < m$. {\rm (}Recall that $\Ga_n$ denotes exactly one of 
$\Ga_n^{\Po}$ or $\Ga_n^{\Q}$.{\rm )} Then there exists some common
$\mu \in \Ga$ such that
$\Ga_{j_i} \mu_{j_i} = \Ga_{j_i} \mu$ for $1 \leq i \leq m$. 
\end{lemma}
\noindent\textbf{Proof}.   Part (a) follows from (\ref{intercon}) just 
as  in \cite[2E]{arp}. 
For part (b) simply note that $\Ga_j \mu = \Ga_k\nu$ forces
$\Ga_j = \Ga_k$ (as with   any two subgroups of a group).

Part (c) is proved just as  in \cite[Lemma 2.2]{harAll}. 
The cases $m = 1,2$ are trivial (even when $j_1 < j_2 = n$),
so assume $m\geq 3$. By induction we have $\gamma$ such that
$\Ga_{j_i} \mu_{j_i} = \Ga_{j_i} \gamma$ for $2 \leq i \leq m$, as well 
as some $\lambda$ such that  
$\Ga_{j_i} \mu_{j_i} = \Ga_{j_i} \lambda$ for $i = 1,2$.
From the overlap at $j_2$ we have $\lambda = \tau \gamma$ for 
some $\tau \in \Ga_{j_2}$. But $j_2 < j_3 \leq n$, so that removal of the 
node labelled $\alpha_{j_2}$ must disconnect the graph in (\ref{2kgp}).
Thus $\tau = \tau' \tau''= \tau'' \tau'$, where $\tau' \in \Ga_{j_i}$ for
all $i = 2, \ldots, m$, and $\tau'' \in \Ga_{j_1} \cap \Ga_{j_2}$.
Then $\mu = (\tau'')^{-1} \lambda=  \tau'\gamma$ is the  desired common coset representative.
\hfill$\square$

 \medskip

\begin{lemma}\label{Spos}
$[\mathcal{S}, \leq]$ is a  partially ordered set
with rank function
\begin{eqnarray*}\rk :  \mathcal{S}  &\rightarrow &\{-1, \ldots, n+1 \}\\
\;\;\;\;\; \; \Ga_j \mu & \mapsto & j
\end{eqnarray*}
\end{lemma}
\noindent\textbf{Proof}.   Keeping Lemma~\ref{dist} in mind, 
it is clear that we need only show that `$\leq$' is transitive.
But that follows at once from the common coset representative 
$\mu$ constructed for part (c) of   Lemma~\ref{dist}. 
\hfill$\square$

\medskip

By Lemma~\ref{dist}(c),
every flag of $\mathcal{S}$ can be written as
$$ [\Ga_0 \mu, \Ga_1\mu, \ldots, \Ga_{n-1}\mu, \Ga_n\mu ]\; ,$$
for some $\mu \in \Ga$ (suppressing the two improper faces). 
Again we emphasize that $\Ga_n$ refers
to exactly one of $\Ga_n^{\Po}$ or $ \Ga_n^{\Q}$, so that there are two 
competing  \textit{base  flags} 
$$ \Phi^{\Po} := [\Ga_0,  \ldots, \Ga_{n-1},   \Ga_n^{\Po}]\;\;\mathrm{and}\;\;
\Phi^{\Q} := [\Ga_0,  \ldots, \Ga_{n-1},  \Ga_n^{\Q}]\; .$$
Clearly there is a right action of $\Ga$ on the flag set $\F(\mathcal{S})$,
and each flag $\Psi$ is equivalent to one  of
$\Phi^{\Po}$ or $\Phi^{\Q}$ 
under this action; on the other hand,
$\Phi^{\Po}$ cannot be $\Ga$-equivalent to $\Phi^{\Q}$, by
Lemma~\ref{dist}(b).
(We later investigate when we can merge these two flag orbits  into one 
by extending $\Ga$.)

Now suppose $K= \{j_1, \ldots, j_m\} $ is some set of ranks satisfying
$0\leq  j_1 < \ldots < j_m\leq n$.  
Each flag $\Psi = \Phi^{\Po} \mu$ contains a chain
$\Psi_K := [\Ga_{j_1} \mu, \ldots, \Ga_{j_m}\mu]$ of type $K$
(with $\Ga_{j_m} = \Ga_n^{\Po}$ if $j_m = n$). Likewise, 
each flag $\Psi = \Phi^{\Q} \mu$ contains a suitable chain $\Psi_K$. 

Notice that the $\Ga$-stabilizer of 
a basic chain $[\Ga_{j_1}, \ldots, \Ga_{j_m}]$ is simply 
\[ \Ga_{K}:=\Ga_{j_1} \cap \ldots \cap \Ga_{j_m}.\]
 By the intersection condition, this is a standard subgroup of $\Gamma$. 
However, a detailed description
of its generators is a bit more involved than for regular polytopes
\cite[Lemma~2E9]{arp}. In particular, as hinted  earlier, here $\Ga_K$ is not simply the subgroup of $\Ga$ generated by all generators but those with $j\not\in K$. 
With this notation in place we have

\begin{lemma}\label{actlemm}
{\rm (a)} The group $\Ga$ acts on  $\F(\mathcal{S})$
and the $\Ga$-stabilizer of the chain  $\Psi_K$ is $\mu^{-1} \Ga_K \mu$.

{\rm (b)} $\Ga$ is transitive on all chains of type $K$ if $n \not\in K$, 
and is transitive on all chains of type $K$ contained in flags 
in either of the two flag orbits.

{\rm (c)} If $K\supseteq \{0,\ldots, n-1\}$, then $\Ga_K$ is trivial.

{\rm (d)} The action of\,  $\Ga$ on $\F(\mathcal{S})$ is faithful,
and we may consider $\Ga$ to be a subgroup of $\Gamma(\mathcal{S})$.
\end{lemma}
\noindent\textbf{Proof}. Part (a) follows easily from our earlier observations.
Part (b) follows at once from Lemma~\ref{dist}(c).
For part (c) we check that $\Ga_0 \cap \ldots \cap \Ga_{n-1}$ is trivial.
In part (d) we need only observe that each $\gamma \in \Ga$ does indeed 
induce an order preserving
bijection on $\mathcal{S}$.
\hfill$\square$

\medskip

\begin{lemma}\label{sects}
{\rm (a)} Let $F$ be any facet of $\mathcal{S}$. Then
the section $F/\Ga_{-1}$ is isomorphic to $\Po$ or to $\Q$.

{\rm (b)} Let $n\geq 2$,  and let $F$ be any vertex of $\mathcal{S}$. Then 
the vertex-figure $\Ga_{n+1}/F$ is isomorphic to  
the ranked poset $\hat{\mathcal{S}}$ constructed in like manner
by deleting the left node in~\emph{(\ref{2kgp})} and transferring the ring to the 
node labelled $\alpha_1$:
\begin{equation}
\label{2kvfig}
\centering
\begin{picture}(180,30)
\put(100,0){
\multiput(15,0)(45,0){2}{\circle*{5}}
\multiput(100,-27)(0,53.3){2}{\circle*{5}}
\multiput(-130,0)(45,0){2}{\circle*{5}}
\put(-130,0){\circle{9}}
\put(-5,0){\line(1,0){20}}
\put(-85,0){\line(1,0){20}}
\put(-130,0){\line(1,0){45}}
\put(15,0){\line(1,0){45}}
\put(60,0){\line(3,2){40}}
\put(60,0){\line(3,-2){40}}
\put(100,-27){\line(0,1){53.3}}
\put(4,9){\scriptsize $\alpha_{n-3}$}
\put(48,9){\scriptsize $\alpha_{n-2}$}
\put(-134,9){\scriptsize $\alpha_{1}$}
\put(-89,9){\scriptsize $\alpha_{2}$}
\put(-49,-0.5){$\ldots\ldots$}
\put(105,-29){\scriptsize $\beta_{n-1}$}
\put(105,25){\scriptsize $\alpha_{n-1}$}
\put(103.5,-3){\scriptsize k}}
\end{picture}
\end{equation}
\vskip.35in

\noindent
{\rm (c) (The base case) } When $n=1$ the  diagram for $\Ga$  becomes  
\begin{equation}
\label{2kpolygon}
\centering
\begin{picture}(180,32)
\put(65,0){
\multiput(0,-27)(0,53.3){2}{\circle*{5}}
\multiput(0,-27)(0,53.3){2}{\circle{9}}
\put(0,-27){\line(0,1){53.3}}
\put(7.5,-29){\scriptsize $\beta_{0}$}
\put(7.5,25){\scriptsize $\alpha_{0}$}
\put(3.5,-3){\scriptsize k}}
\end{picture}
\end{equation}
\vskip.35in
\noindent
which describes the polygon $\mathcal{S}=\{2k\}$. 
\end{lemma}

\noindent\textbf{Proof}.  In part (a) we may assume
without loss of generality  
that $F = \Ga_n^{\Po}$. Since $\Ga_{n-1} \subset \Ga_n^{\Po}$, an $(n-1)$-face
$\Ga_{n-1}\mu \leq \Ga_n^{\Po}$ if and only if 
$\mu \in \Ga_n^{\Po}$. By Lemma~\ref{dist}(c), this means that all $j$-faces
$\Ga_j\mu$ incident with 
$\Ga_n^{\Po}$ are likewise represented by certain $\mu \in \Ga_n^{\Po}$.
On the other hand, a basic $j$-face of $\Po$ can be identified with 
\[ G_j := \langle \alpha_0, \ldots, \alpha_{j-1}, \alpha_{j+1},
\ldots, \alpha_{n-1}\rangle = \Ga_j \cap \Ga_n^{\Po}.\] 
We therefore have a well-defined bijection 
$$\begin{array}{rcl}
\eta: \Po\;\;\; & \rightarrow & \Ga_n^{\Po}/\Ga_{-1}\\
G_j \mu & \mapsto & \Ga_j \mu\; ,
\end{array}$$
for  $\mu \in \Ga_n^{\Po}$, $0 \leq j \leq n-1$. It is easy   to check that 
$\eta$  is  order-preserving. Now suppose that 
$\Ga_j\mu  \cap \Ga_k\nu \neq \emptyset$,
for $j <k$ and $\mu, \nu \in \Ga_n^{\Po}$. By Lemma~\ref{dist}(c) we may assume 
$\mu = \nu$;  clearly $G_j \mu \cap G_k\mu \neq \emptyset$.
Thus $\eta^{-1}$ is also order-preserving.

Now we turn to part (b). We may assume that $F = \Ga_0$. Notice that
$\Ga_0$ is nothing more than the group described by the diagram
(\ref{2kvfig}) and that this group still satisfies the 
intersection condition on its own standard subgroups.
Thus the proper $k$-faces of the ranked poset $\hat{\mathcal{S}}$ are 
cosets of suitable basic subgroups. In fact, by 
the intersection condition for $\Ga$
we find that  these basic subgroups are
$$\Ga_{0,k}:= \Ga_0 \cap \Ga_{k+1} \; ,$$
for $0 \leq k \leq n-1$. (Thus $\Ga_{0,n-1}$ comes in 
two varieties, as determined by
the vertex-figures of $\Po$ and $\Q$.)

Once more using Lemma~\ref{dist}(c) we find for $j \geq 1$ that any $j$-face $\Ga_j\mu$ incident with $\Ga_0$  
is represented by some $\mu \in \Ga_0$. We may therefore define
\begin{eqnarray*}
\lambda : \Ga_{n+1}/\Ga_0 & \rightarrow & \hat{\mathcal{S}}\\
\Ga_j \mu & \mapsto & \Ga_{0,j-1}\,  \mu = (\Ga_0 \cap \Ga_j)\mu\; ,
\end{eqnarray*}
for $\mu \in \Ga_0$, $1\leq j \leq n$ (again ignoring improper faces).
As before, it is easy to check that $\lambda$ is a poset isomorphism.

In part (c), the group $\Ga = \langle \alpha_0, \beta_0 \rangle$ is the dihedral 
group $\mathbb{D}_{2k}$.                  
There are two basic facets (here edges), namely 
$\Ga_1^{\Po} = \langle \alpha_0 \rangle$ and 
$\Ga_1^{\Q} = \langle \beta_0\rangle$, each incident with the basic
vertex $\Ga_0 = \{1\}$. One checks that these faces are arranged as indicated 
in (\ref{2sect}) below (taking $n=1$). 
\hfill$\square$

\begin{theorem} \label{wytsemi}
The abstract Wythoff's construction described in  Definition~\ref{defS}
and summarized in diagram
{\rm (\ref{2kgp})}  defines a semiregular $(n+1)$-polytope $\mathcal{S}$.
Its facets are isomorphic to $\Po$ or $\Q$, with $k$ of each of these occurring 
alternately around each face
$\R$ of co-rank $2$. {\rm (}Thus each $2$-section 
$\mathcal{S}/\R$ is a $2k$-gon.{\rm )}
The face-wise  $\Ga$-stabilizer of any ridge of $\mathcal{S}$ is trivial. 
Finally, each vertex-figure of $\mathcal{S}$
is isomorphic to  the semiregular $n$-polytope $\hat{\mathcal{S}}$ defined 
by the diagram
in \emph{(\ref{2kvfig})}.
\end{theorem}

\noindent\textbf{Proof}. We adapt the methods  of \cite[2E]{arp}.
Polytope property \textbf{A}  is already  in place, so we move to the 
`diamond condition' \textbf{B}. 
For $0 \leq j \leq n$, 
we must consider incident faces
of ranks $j-1$ and $j+1$. By  Lemma~\ref{actlemm}(b)   we can 
take these to be $\Ga_{j-1}$ and $\Ga_{j+1}$ in the basic chain
$$ \Phi_K = [\Ga_{-1}, \ldots, ,\Ga_{j-1}, \Ga_{j+1}, \ldots,\Ga_n,\Ga_{n+1}]\;, 
$$
with $K:=\{0,1,\ldots,n+1\}\setminus \{j\}$.  

If $j\leq n-2$, the stabilizer of
this chain is given by $\Ga_K = \langle \alpha_{j} \rangle = \{1,\alpha_j\}$.
This implies that the distinct faces $\Ga_{j}$ and $\Ga_{j}\alpha_j$
are the only $j$-faces incident with $\Ga_{j-1}$ and $\Ga_{j+1}$.

If $j = n-1$, then $\Ga_K = \langle \alpha_{n-1} \rangle$ (resp.
$\langle \beta_{n-1} \rangle$) if $\Ga_n = \Ga_n^{\Po}$ 
(resp. $\Ga_n = \Ga_n^{\Q}$ ), and we proceed similarly.

If $j = n$, then $\Ga_K = \{1\}$. However, 
$\Ga_{n-1} = \Ga_n^{\Po} \cap \Ga_n^{\Q}$, so 
$\Ga_{n-1} \cap \Ga_n^{\Po}\mu \neq\emptyset$ forces 
$\Ga_n^{\Po}\mu  = \Ga_n^{\Po}$. Thus $\Ga_n^{\Po}$ and $\Ga_n^{\Q}$
are the two distinct facets incident with the ridge $\Ga_{n-1}$.

Next let us examine the structure of the $2$-section $\mathcal{S}/\R$. 
Without loss of generality
we can assume that $\R = \Ga_{n-2}$.  It is easy to check that 
the $(n-1)$-faces in this section are   all $\Ga_{n-1}\mu$, 
with $\mu$ in the dihedral group 
$\mathbb{D}_{2k} := \langle \alpha_{n-1}, \beta_{n-1} \rangle$.
Since $\mathbb{D}_{2k}  \cap \Ga_{n-1} = \{1\}$ (by the intersection condition),
there are exactly $2k$ such cosets.
Similarly, the $n$-faces come in two types:
$k$ each of types  $\Ga_n^{\Po}\mu$ and 
$\Ga_n^{\Q}\mu$, again taking  $\mu \in  \mathbb{D}_{2k}$. 
Thus each $\Ga_{n-1}\mu$ meets exactly one $n$-face of each type,
namely $\Ga_n^{\Po}\mu$ and $\Ga_n^{\Q}\mu$. Likewise,
$\Ga_n^{\Po}\mu$ (resp. $\Ga_n^{\Q}\mu$) meets just
 $\Ga_{n-1}\mu$ and  $\Ga_{n-1}\alpha_{n-1}\mu$
(resp.  $\Ga_{n-1}\mu$ and  $\Ga_{n-1}\beta_{n-1}\mu$).
We see that this section has the structure of a $2k$-gon, as indicated in 

\begin{equation}
\label{2sect}
\centering
\begin{picture}(400,12)
\put(-1,-12){
\put(-8,14.5){$\ldots$}
\put(401,14.5){$\ldots$}
\multiput(75,15)(65,0){5}{\circle*{5}}
\multiput(10,15)(130,0){3}{\line(1,0){65}}
\put(75,15){
\multiput(5,0)(10,0){6}{\line(1,0){5}}
\multiput(135,0)(10,0){6}{\line(1,0){5}}
\multiput(265,0)(10,0){6}{\line(1,0){5}}}
\put(50,3){\scriptsize $\Ga_{n-1}\beta_{n-1}\alpha_{n-1}$}
\put(125,3){\scriptsize $\Ga_{n-1}\alpha_{n-1}$}
\put(195,3){\scriptsize $\Ga_{n-1}$}
\put(255,3){\scriptsize $\Ga_{n-1}\beta_{n-1}$}
\put(310,3){\scriptsize $\Ga_{n-1}\alpha_{n-1}\beta_{n-1}$}
\put(10,22){\scriptsize $\Ga_{n}^{\Po}\beta_{n-1}\alpha_{n-1}$}
\put(93,22){\scriptsize $\Ga_{n}^{\Q}\alpha_{n-1}$}
\put(166,22){\scriptsize $\Ga_{n}^{\Po}$}
\put(231,22){\scriptsize $\Ga_{n}^{\Q}$}
\put(288,22){\scriptsize $\Ga_{n}^{\Po}\beta_{n-1}$}
\put(343,22){\scriptsize $\Ga_{n}^{\Q}\alpha_{n-1}\beta_{n-1}$}}
\end{picture}
\end{equation}
\vskip.2in

Now we can verify the strong connectedness property \textbf{C}.   
We must show that any section 
$F/G$ of $\mathcal{S}$ is connected. If $\rk(G) \geq 0$ or $\rk(F) \leq n$, 
this follows from a standard inductive argument based on Lemma~\ref{sects}
and the fact that $\Po$ and $\Q$, being   polytopes, are 
themselves strongly connected. To complete the induction we need only show that
$\mathcal{S}$ is itself  connected. For this it will suffice 
to show that any facet $\Gamma_n^{\Po}\mu$ can be connected to 
$\Gamma_n^{\Po}$, say, 
by a sequence of consecutively incident facets and ridges.
If $\mu = \beta_{n-1}$, we observe such a sequence in (\ref{2sect}) above.
(Of course, facets of type $\Ga_n^{\Q}$ will appear in the sequence.)
Translating by any $\gamma \in \Ga$ we get a sequence connecting 
$\Ga_n^{\Po} \gamma$ to $\Ga_n^{\Po}\beta_{n-1} \gamma$.

Now consider a general facet $\Ga_n^{\Po} \mu$, and write
$\mu$ as a word of minimal length in the generators, say
$$ \mu = \rho_1 \cdots \rho_m, \;\;\mathrm{with}\;\; \rho_j \in 
\{ \alpha_0, \ldots, \alpha_{n-1},\beta_{n-1} \}\; .$$
Suppose $m=1$. If $\mu = \rho_1 = \alpha_j$ for some $j$, then 
$\Ga_n^{\Po}\mu = \Ga_n^{\Po}\alpha_j =   \Ga_n^{\Po}$ and there is nothing to
prove. If $\mu = \rho_1 = \beta_{n-1}$, then we have already observed the
required sequence in (\ref{2sect}) above.

We proceed inductively, assuming the existence of a suitable sequence 
connecting $\Ga_n^{\Po}$ to  $\Ga_n^{\Po}\gamma$ whenever
$\gamma$ has length $m-1$ in the generators.
Let $\mu = \rho_1 \cdots \rho_m$, so that $\mu =\rho_1 \gamma$ for
$\gamma = \rho_2 \cdots \rho_m.$ If $\rho_1 = \alpha_j$, then as before
$\Ga_n^{\Po}\mu = \Ga_n^{\Po}\gamma$, which by induction can 
in fact be joined to $\Ga_n^{\Po}$.
On the other hand, if $\rho_1 = \beta_{n-1}$,
then 
$\Ga_n^{\Po}\mu = \Ga_n^{\Po}\beta_{n-1} \gamma$
is linked as in (\ref{2sect})  to $\Ga_n^{\Po}\gamma$,
which in turn is linked to 
$\Ga_n^{\Po}$ by the induction hypothesis.

This completes the proof.
\hfill$\square$

\medskip

\noindent\textbf{Remarks}. 
In the particularly interesting special issue when the diagram in (\ref{2kgp}) 
is a star with three branches (that is $n=3$ and $k=2$),  there generally are 
three non-equivalent ways in which the diagram yields semiregular polytopes of 
rank $4$. For example, the affine Coxeter group $\widetilde{B}_3$ represented 
by the diagram in (\ref{fig2}) actually gives rise to two semiregular $3$-polytopes, 
namely the semiregular tessellation of $\mathbb{R}^3$ with tetrahedral and 
octahedral tiles shown in Figure~\ref{fig1}, as well as the standard cubical 
tessellation of $\mathbb{R}^3$; the latter is geometrically regular, by our 
next proposition, and $\widetilde{B}_3$ is a subgroup of index $2$ in its full
 symmetry group $\widetilde{C}_3$. Here, the cubical tessellation can be 
derived in two ways from the diagram; on other words, two of the three 
semiregular polytopes associated with the diagram actually are isomorphic. 
In Section~\ref{finfi} we describe examples where the three semiregular 
polytopes are mutually non-isomorphic. At the other extreme, for the finite 
Coxeter group $D_4$,  
 all three semiregular polytopes are mutually isomorphic; 
in fact, each is regular and isomorphic to the $4$-cube.  

Recall that an abstract polytope is called a {\em $2$-orbit polytope\/} if its automorphism 
group has precisely two flag orbits. 
(See \cite{hub2orb}; for  
general structure results about the groups of $2$-orbit polytopes see also \cite{hubsch}.) 

\begin{proposition}
\label{2orbit}
Suppose $\mathcal{S}$ is the semiregular 
$(n+1)$-polytope  constructed from the tail-triangle diagram in \emph{(\ref{2kgp})}.

\emph{(a)} Then $\mathcal{S}$ is a regular  polytope if and only if $\Gamma$ admits a group automorphism induced by the diagram symmetry which swaps $\alpha_{n-1}$ and $\beta_{n-1}$ in \emph{(\ref{2kgp})}, while fixing the remaining $\alpha_j$'s. In this case $\Po \simeq \Q$, say with Schl\"{a}fli type  $\{p_1, \ldots, p_{n-1}\}$, and $\mathcal{S}$ is regular 
of type $\{p_1, \ldots, p_{n-1}, 2k \}$; moreover, $\Gamma(\mathcal{S}) \simeq \Gamma  \rtimes C_2$.

\emph{(b)} If $\mathcal{S}$ is not regular,  then $\mathcal{S}$ is a 2-orbit polytope 
and $\Gamma(\mathcal{S}) \simeq \Ga$. In particular, this is so if the facets $\Po$ and $\Q$ are non-isomorphic (as is the case, for example, if $\alpha_{n-2} \, \alpha_{n-1}$ and  $\alpha_{n-2}\,  \beta_{n-1}$ 
have different periods). 

\end{proposition}
 
\noindent\textbf{Proof}. 
Each flag of $\mathcal{S}$ is equivalent under $\Ga$ to exactly one of the base flags 
\[ \Phi^{\Po} := [\Ga_0,  \ldots, \Ga_{n-1},   \Ga_n^{\Po}], \;\;\mbox{or} \;\;
\Phi^{\Q} := [\Ga_0,  \ldots, \Ga_{n-1},  \Ga_n^{\Q}]\; \]
(relative to $\Ga$). Thus $\mathcal{S}$ has two flag orbits under $\Ga$, and one or two 
flag orbits under 
its full automorphism group $\Gamma(\mathcal{S})$. In particular, 
$\mathcal{S}$ is regular if and only if $\Gamma(\mathcal{S})$ has just one flag orbit, 
and then $\Ga$ has index $2$ in $\Gamma(\mathcal{S})$. 
%

Now suppose $\mathcal{S}$ is regular and $\Phi:=\Phi^{\Po}$ denotes the base flag of 
$\mathcal{S}$ relative to the full group $\Gamma(\mathcal{S})$. Then 
$\Gamma(\mathcal{S}):=\langle \gamma_0,\ldots,\gamma_{n}\rangle$, where 
$\gamma_0,\ldots,\gamma_{n}$ are the distinguished generators. Identifying 
$\Ga$ with a subgroup of $\Gamma(\mathcal{S})$ and inspecting the  action of its 
generators on $\Phi$, we find that $\gamma_{j}=\alpha_{j}$ for $j\leq n-1$. Hence, 
if $j\leq n-2$ then 
\[ \gamma_{n}\alpha_{j}\gamma_{n} = \gamma_{n}\gamma_{j}\gamma_{n} = \gamma_{j}=\alpha_{j}. \]
Moreover, conjugation in $\Gamma(\mathcal{S})$ by $\gamma_n$ takes the distinguished 
generators $\gamma_0,\ldots,\gamma_{n}$ relative to $\Phi=\Phi^{\Po}$ to the 
distinguished generators relative to the adjacent flag $\Phi^{n}=\Phi^{Q}$ of $\mathcal{S}$. 
Arguing as before we see from the action on $\Phi^{\Q}$ that 
$\gamma_{n}\gamma_{n-1}\gamma_{n-1}=\beta_{n-1}$ and therefore
\[ \gamma_{n}\alpha_{n-1}\gamma_{n} = \gamma_{n}\gamma_{n-1}\gamma_{n} = \beta_{n-1}.\]
Hence conjugation in $\Gamma(\mathcal{S})$ by $\gamma_n$ induces an involutory group 
automorphism on $\Ga$ corresponding to the horizontal diagram symmetry. Clearly, 
$\Gamma(\mathcal{S}) \simeq\Gamma  \rtimes C_2$, with $C_{2} = \langle \gamma_n\rangle$. 
The remaining claims of part (a) are obvious.

Conversely, suppose that an automorphism $\tau$ of $\Gamma$ is induced by the diagram 
symmetry which swaps $\alpha_{n-1}$ and $\beta_{n-1}$, while fixing the remaining 
$\alpha_j$'s. Then, in the terminology of \cite[8A]{arp}, a twisting operation applied 
with $\Gamma$ and $\tau$ recovers the polytope $\mathcal{S}$ and proves that it is regular. 
The details are routine.

Now suppose $\mathcal{S}$ is not regular (for example, this occurs if the two facet types 
$\Po$ and $\Q$ are not isomorphic). Then both $\Gamma(\mathcal{S})$ and its subgroup $\Ga$ 
have two flag orbits, so $\Gamma(\mathcal{S}) \simeq \Ga$ and $\mathcal{S}$ is a $2$-orbit 
$(n+1)$-polytope (of type $2_{\{0, \ldots,n-1\}}$, in the terminology of~\cite{hub2orb}). 
\hfill$\square$
\medskip

\noindent\textbf{Remarks}. 
We note for part (b) that $\Ga =  \Gamma(\mathcal{S})$ can certainly
hold even when $\Po \simeq \Q$. There are already instances of this 
behaviour when $n = 2$, for we need only  truncate a regular 
polyhedron of type $\{p,p\}$
which is \textit{not} self-dual. A quick check of Hartley's 
Census \cite{harcen} reveals several instances.
The smallest such polyhedron is 
flat with $72$ flags and 
Schl\"{a}fli type $\{6,6\}$. However,
for ease of description we choose instead the following 

\begin{example}\label{66map}\emph{ 
Let $\mathcal{M}$ to be the regular polyhedron
$\{6,6\}\ast 240a$  in the Census \cite{harcen}.  
The corresponding string C-group 
$\Gamma$ has
order $240$ and can be
generated by permutations
$$\alpha_1 = (2,3)(4,5), \; \alpha_0 = (1,2),\; \beta_1 = (2,4)(3,5)(6,7)$$
associated to the diagram
\begin{equation}
\centering
\begin{picture}(180,30)
\put(0,0){
\multiput(100,-27)(0,53.3){2}{\circle*{5}}
\put(60,0){\circle*{5}}
\put(60,0){\circle{9}}
\put(60,0){\line(3,2){40}}
\put(60,0){\line(3,-2){40}}
\put(77,16){\scriptsize 6}
\put(77,-21.5){\scriptsize 6}
\put(43,-2){\scriptsize $\alpha_{0}$}
\put(105,-29){\scriptsize $\beta_{1}$}
\put(105,25){\scriptsize $\alpha_{1}$}}
\end{picture}
\end{equation}
\vskip.4in
\noindent
We construct $\mathcal{M}$  by ringing the top node and
its  dual  $\mathcal{M}^\ast$ by ringing the bottom node.
If, however, we  ring the middle node, then  we obtain a semiregular 
polyhedron $\mathcal{S}$ with $480$ flags. This $\mathcal{S}$
has Schl\"{a}fli type  $\{6,6\}$, but it cannot be regular and
$\Gamma(\mathcal{S}) \simeq \Ga$ has order $240$.
Intuitively, we construct $\mathcal{S}$   by truncating 
$\mathcal{M}$  to its edge midpoints.
}\end{example}

\begin{example} \emph{The Coxeter group
$B_3$ of order $48$ has diagram
\begin{equation}
\centering
\begin{picture}(180,30)
\put(0,0){
\multiput(100,-27)(0,53.3){2}{\circle*{5}}
\put(60,0){\circle*{4}}
\put(60,0){\circle{9}}
\put(60,0){\line(3,2){40}}
\put(100,-27){\line(0,1){53.3}}
\put(77,16){\scriptsize 4}
\put(43,-2){\scriptsize $\alpha_{0}$}
\put(105,-29){\scriptsize $\beta_{1}$}
\put(105,25){\scriptsize $\alpha_{1}$}
\put(102.5,-3){\scriptsize 3}}
\end{picture}
\end{equation}
\vskip.4in 
\noindent 
The geometrical version of Wythoff's construction gives 
the cube $\{4,3\}$, which is a  convex regular solid, of course. 
The abstract construction is  subtly different. Since $\alpha_0 \beta_1$
has period $2$,  each edge of the cube must be replaced by a digon, so that 
we do indeed get a semiregular polytope with  square and 
digonal faces and with hexagonal vertex-figures.
In the convex setting, these digons collapse to line segments and
the vertex-figures to equilateral triangles, as expected.
(A nice spherical model is obtained by inscribing the cube in its 
circumsphere; the six face planes cut the sphere in the
$24$ arcs needed for edges.) }              

\emph{In the same way we can disturb any of the classical 
convex regular 
polytopes, indeed, any regular polytope whatsoever. }  
\end{example}
\medskip

The tessellation $\T$ described in
 Example~\ref{tessT} and Figure~\ref{fig1}
has octahedra  among its facets, and cuboctahedral vertex-figures.
Both have spherical type and can be projectified by antipodal
identifications. This is accomplished by adjoining to (\ref{B3cox})
the first or second relation in (\ref{tomrels}) below. We choose both for the 
next example (see \cite{monpelwillA} for the 
other possibilities):

\begin{example}\label{tomeg} \textbf{The tomotope}.
\emph{
Define a group $\Ga$ by adjoining to (\ref{B3cox})
the two relations which create hemioctahedra and hemicuboctahedra:
}
\begin{equation}\label{tomrels}
 (\rho_0 \rho_1 \rho_3)^3 = (\rho_2 \rho_1 \rho_3)^3 = 1\; .
 \end{equation}
\end{example}
\noindent
Using GAP \cite{gap}, it is easy to check
that $\Ga$ now has finite order $96$, yet still satisfies the intersection
condition (\ref{intercon}). Thus, by Theorem~\ref{wytsemi}, we obtain
a finite  semiregular $4$-polytope $\tom$, which in \cite{monpelwillA} 
was called the \textit{tomotope}
(as a small gift for Toma\v{z} Pisanski).

The tomotope has $4$ vertices, $12$ edges, $16$ triangles, and $4$ 
tetrahedra and $4$ hemioctahedra. Two of each kind of facet
alternate around each edge. Since $\Ga$ acts faithfully on edges
 $1,\ldots,12$ (say),
we have this permutation representation:

\begin{eqnarray*}
\rho_0 &= &(5,10) (6,9)(7,12)(8,11)\,,\\
\rho_1 &= &(1,6) (2,5)(3,8) (4,7)\,,\\ 
\rho_2 &= &(5,9) (6,10) (7, 11) (8,12)\,,\\
\rho_3 &= &(5,8) (6,7) (9,12) (10,11)\, .
\end{eqnarray*}
It is also possible to obtain $\Ga$ by reducing the 
crystallographic group $\widetilde{B}_3$ modulo $2$.
(This fact resurfaces in Section \ref{finfi} below.)

Our main purpose in \cite{monpelwillA}  was to investigate  
the regular 
covers of $\tom$. Since the epimorphism 
$\widetilde{B}_3 \rightarrow \Ga$ induces a $2$-covering 
$ \T \rightarrow \tom$
of semiregular polytopes \cite[2D]{arp},  we may construct $\tom$  
by making suitable identifications in $\T$ (see below).  
 But putting that aside,  we find that
the tomotope has the peculiar property, impossible
in lower ranks, of having infinitely many distinct, finite, minimal
\textit{regular} covers.

\medskip 
 
\begin{figure}[htbp]
\centering
\includegraphics[width=60mm]{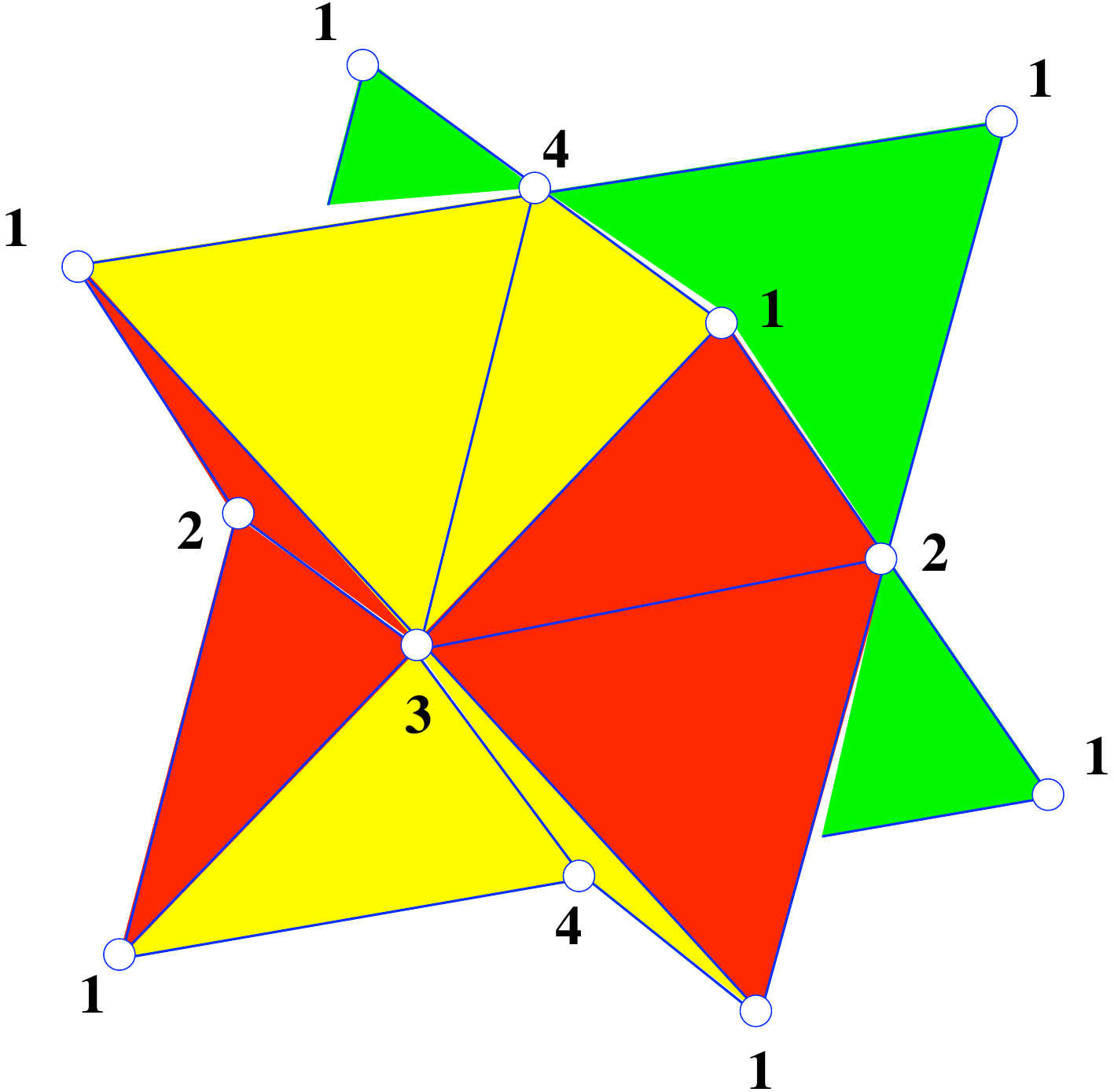}
\caption{The tomotope $\tom$.}\label{fig4a}
\end{figure}

To visualize the tomotope $\tom$
imagine a core octahedron with $8$ tetrahedra
glued to its faces, suggesting the \textit{stella octangula}.
Next  imagine this complex inscribed in a $2 \times 2 \times 2$
cube and from that
make toroidal-type identifications  for the vertices 
and edges lying in the boundary of the cube. Finally, we further 
identify antipodal faces of all ranks  to get $\tom$.
In Figure~\ref{fig4a} you can
see the $4$ vertices,  $4 = 8/2$ tetrahedra and $1$  
hemioctahedron in the core. The other three hemioctahedra are red, 
yellow  and green,
and `run around' the belts of those colours. For example,
before making identifications, we may split a red octahedron into four 
sectors around a vertical axis of symmetry, then fit these into  
four red slots, two of which are visible
in Figure~\ref{fig4a}. In this way we fill out the $2 \times 2 \times 2$
cube before the final antipodal identifications.

The proof of the intersection  property for tail-triangle C-groups can often be 
reduced to the consideration of only a small number of cases. For 
$i=-1,0,\ldots,n-2$,  define 
$\Ga_{i}^+:=\langle\alpha_{i+1},\ldots,\alpha_{n-2},\alpha_{n-1},\beta_{n-1}\rangle$. 

\begin{lemma}
\label{intersec}
Suppose that $\Ga = \langle \alpha_0, \ldots, \alpha_{n-2},\alpha_{n-1},\beta_{n-1}\rangle$ is a tail-triangle group with $n\geq 2$, and suppose that its subgroups $\Ga_n^1:=\langle \alpha_0, \ldots, \alpha_{n-2},\alpha_{n-1}\rangle$, $\Ga_n^2:=\langle \alpha_0, \ldots, \alpha_{n-2},\beta_{n-1}\rangle$ and $\Ga_{0}:=\langle \alpha_1, \ldots, \alpha_{n-2},\alpha_{n-1},\beta_{n-1}\rangle$ are C-groups. Then $\Ga$ is a tail-triangle C-group if and only if $\Ga_{n}^{1}\cap\Ga_{n}^{2} = \langle \alpha_0, \ldots, \alpha_{n-2}\rangle$ and both $\Ga_{i}^{+}\cap\Ga_{n}^1 = \langle \alpha_{i+1}, \ldots, \alpha_{n-1}\rangle$ and $\Ga_{i}^{+}\cap\Ga_{n}^2 = \langle \alpha_{i+1}, \ldots, \alpha_{n-2},\beta_{n-1}\rangle$ for $i=0,\ldots,n-2$.
\end{lemma}

\noindent\textbf{Proof}. 
The case $n=2$ is trivial,  since the stated conditions on 
$\Ga_{n}^1$, $\Ga_{n}^2$ and $\Ga_{0}^+$ are just saying that any 
two of the three distinguished dihedral subgroups intersect as required. 
Now suppose $n\geq 3$ and 
$I,J\subseteq \{\alpha_0,\ldots,\alpha_{n-1}, \beta_{n-1}\}$. 
We must show that 
$\langle I \rangle \cap \langle J \rangle= \langle I \cap J\rangle$. 
Clearly, one inclusion is obvious, and the other holds trivially if 
$I$ or $J$ is the full generating set. Hence we may exclude the latter 
possibility from now on. 

First consider the case when $\{\alpha_0,\ldots,\alpha_{n-2}\}\not\subseteq I,J$. 
Let $i$ (resp. $j$) denote the largest integer $k$ with $k\leq n-2$ and 
$\alpha_{k}\not\in I$ (resp. $\alpha_{k}\not\in J$). Assume that $i\leq j$. 
Then $I = I_0 \cup I_1$ and $J = J_0 \cup J_1$, where 
\[ I_0\subseteq \{\alpha_0, \ldots,\alpha_{i-1}\},\;\;\; 
\{\alpha_{i+1},\ldots,\alpha_{n-2}\} \subseteq I_1\subseteq \{\alpha_{i+1},\ldots,\alpha_{n-1}, \beta_{n-1}\} \] 
and similarly 
\[ J_0\subseteq \{\alpha_0, \ldots,\alpha_{j-1}\},\;\;\; 
\{\alpha_{j+1},\ldots,\alpha_{n-2}\} \subseteq J_1\subseteq \{\alpha_{j+1},\ldots,\alpha_{n-1}, \beta_{n-1}\}.\] 
By the commutativity relations implicit in the underlying diagram, 
$\langle I \rangle \simeq \langle I_0 \rangle\times \langle I_1 \rangle$ and 
 $\langle J \rangle \simeq \langle J_0 \rangle\times \langle J_1 \rangle$.
(It will emerge below that these and similar products are direct.) 

Now suppose $\mu\in \langle I \rangle \cap \langle J \rangle$. Then 
$\mu = \lambda_{0}\lambda_{1}$ with $\lambda_{0}\in \langle I_0\rangle$ 
and $\lambda_{1}\in \langle I_1\rangle$, and similarly 
$\mu = \nu_{0}\nu_{1}$ with $\nu_{0}\in \langle J_0\rangle$ and 
$\nu_{1}\in \langle J_1\rangle$. But $i\leq j$, so both $\lambda_{1}$ and $\nu_1$ are 
elements in $\langle \alpha_{i+1},\ldots,\alpha_{n-1}, \beta_{n-1} \rangle$ related by 
$\nu_{1}=\xi \lambda_{1}$ with 
$\xi:= \nu_{0}^{-1}\lambda_{0}\in \langle \alpha_0, \ldots,\alpha_{j-1} \rangle$ 
(and $\xi=\nu_{1}\lambda_{1}^{-1} \in\langle \alpha_{i+1},\ldots,\alpha_{n-1}, 
\beta_{n-1}\rangle$). Then 
by our assumptions 
$\xi\in \Ga_{i}^{+}\cap\Ga_{n}^1 = \langle \alpha_{i+1}, \ldots,\alpha_{n-1}\rangle$, 
and hence also 
$\xi\in \langle \alpha_0, \ldots,\alpha_{j-1} \rangle \cap \langle \alpha_{i+1}, 
\ldots,\alpha_{n-1}\rangle
= \langle \alpha_{i+1}, \ldots,\alpha_{j-1}\rangle$ since $\Ga_{n}^1$ is a C-group. 
(Replacing $I,J$ by $I_0, I_1$ we now see why 
$\langle I_0\rangle \cap \langle I_1\rangle = \langle 1\rangle$.) 
Returning to our proof, we note that
it suffices to show that $\lambda_0,\lambda_1\in 
\langle I \cap J \rangle$.

First consider $\nu_0=\lambda_0\xi^{-1}$ as an element of the direct product 
\[\langle \alpha_{0},\ldots,\alpha_{i-1},\alpha_{i+1},\ldots,\alpha_{j-1}\rangle \simeq
\langle \alpha_{0},\ldots,\alpha_{i-1}\rangle \times \langle\alpha_{i+1},\ldots,
\alpha_{j-1}\rangle,\]
bearing in mind here that $\lambda_0\in\langle \alpha_{0},\ldots,\alpha_{i-1}\rangle $ 
and $\xi\in \langle\alpha_{i+1},\ldots,\alpha_{j-1}\rangle$. But 
\[ \nu_{0}\in \langle J_0\rangle \cap \langle \alpha_{0},\ldots,\alpha_{i-1},
\alpha_{i+1},\ldots,\alpha_{j-1}\rangle      
=\langle J_0\setminus\{\alpha_i\}\rangle,\] 
once again 
by the intersection property in $\Ga_{n}^1$. Since  the 
factorization of the element $\nu_0$ of $\langle J_0\setminus\{\alpha_i\}\rangle$ as an 
element of the direct product is unique,  
we must also have $\lambda_0,\xi\in \langle J_0\rangle$; in fact, we even have 
$\xi\in \langle J_0\setminus \{\alpha_0\}\rangle$. But then $\lambda_0\in \langle 
I_0\rangle \cap \langle J_0\rangle = \langle I_0 \cap J_0\rangle \leq \langle I 
\cap J\rangle$, by the intersection property in $\Ga_{n}^1$.

Similarly we look at the factorization $\lambda_1 = \xi^{-1}\nu_1$ in the direct product 
\[\langle \alpha_{i+1},\ldots,\alpha_{j-1},\alpha_{j+1},\ldots,\alpha_{n-1},\beta_{n-1}
\rangle \simeq \langle \alpha_{i+1},\ldots,\alpha_{j-1}\rangle \times \langle\alpha_{j+1},
\ldots,\alpha_{n-1},\beta_{n-1}\rangle,\]
remembering here that $\xi\in \langle\alpha_{i+1},\ldots,\alpha_{j-1}\rangle$ and 
$\nu_1\in\Ga_{j+1}^{+}$. Now   
\[ \lambda_{1}\in \langle I_1 \rangle \cap 
\langle \alpha_{i+1},\ldots,\alpha_{j-1},\alpha_{j+1},\ldots,\alpha_{n-1},\beta_{n-1}
\rangle
= \langle I_1\setminus \{\alpha_j\}\rangle, \]
by the intersection property of $\Ga_0$. But the factorization of the element 
$\lambda_{1}$ of $\langle I_1\setminus \{\alpha_j\}\rangle$ 
as an element of the direct 
product is unique, so in fact $\xi,\nu_1\in \langle I_1 \rangle$. Hence,   
\[\nu_1\in \langle I_1 \rangle \cap  \langle J_1 \rangle =  \langle I_1 \cap J_1 
\rangle \leq \langle I \cap J\rangle,\] 
by the intersection property of $\Ga_0$. Moreover, 
\[\xi\in \langle I_1 \rangle \cap \langle J_0\setminus \{\alpha_0\}\rangle = 
\langle I_1 \cap (J_0\setminus \{\alpha_0\})\rangle\leq \langle I \cap J\rangle,\] 
once again by the intersection property of $\Ga_0$. 
Thus $\lambda_1\in \langle I \cap J\rangle$, as required. This completes the proof for 
the case when $\{\alpha_0,\ldots,\alpha_{n-2}\}\not\subseteq I,J$.

Now suppose that one of the sets, say $J$, 
contains $\{\alpha_0,\ldots,\alpha_{n-2}\}$. 
Then either $J\subseteq \{\alpha_0,\ldots,\alpha_{n-1}\}$ and 
$\langle J \rangle \leq \Ga_{n}^1$, 
or $J=\{\alpha_0,\ldots,\alpha_{n-2},\beta_{n-1}\}$ and $\langle J \rangle = \Ga_{n}^2$. 

If also $I$ contains $\{\alpha_0,\ldots,\alpha_{n-2}\}$, then either $I\subseteq J$ or 
$J\subseteq I$, or one of $I,J$ is $\{\alpha_0,\ldots,\alpha_{n-1}\}$ and the other is 
$\{\alpha_0,\ldots,\alpha_{n-2},\beta_{n-1}\}$. In the former case the intersection 
condition holds trivially, and in the latter case by our assumptions on $\Ga_n^1 \cap 
\Ga_{n}^2$. 

Now suppose that $\{\alpha_0,\ldots,\alpha_{n-2}\}\not\subseteq I$ and that $\mu\in \langle I \rangle \cap \langle J \rangle$. As above let $i$ denote the largest integer $k$ with $k\leq n-2$ and $\alpha_{k}\not\in I$, let $I = I_0 \cup I_1$, and let $\mu = \lambda_{0}\lambda_{1}$ with $\lambda_{0}\in \langle I_0\rangle$ and $\lambda_{1}\in \langle I_1\rangle$. 
But $I_0\subseteq J$, so clearly $\lambda_0 \in\langle I \cap J \rangle$ and $\lambda_{1}=\lambda_0^{-1}\mu\in \langle J\rangle$.

First suppose that $J\subseteq \{\alpha_0,\ldots,\alpha_{n-1}\}$. Then 
\[ \lambda_{1}\in \langle I_1\rangle \cap \langle J\rangle \leq \Ga_{i}^{+} \cap\Ga_n^1
= \langle\alpha_{i+1},\ldots,\alpha_{n-1}\rangle, \]
by our assumption on the rightmost intersection. Hence, if $\alpha_{n-1}\in I_1,J$ then  necessarily $\{\alpha_{i+1},\ldots,\alpha_{n-1}\} = I_{1}\cap J$ and hence $\lambda_{1}\in \langle I\cap J \rangle$; bear in mind that 
$\{\alpha_{i+1},\ldots,\alpha_{n-2}\} \subseteq  I_1$ by the definition of $i$. Now, if $\alpha_{n-1}\not\in J$ then 
\[ \lambda_{1}\in \langle J \rangle \cap \langle\alpha_{i+1},\ldots,\alpha_{n-1}\rangle 
= \langle\alpha_{i+1},\ldots,\alpha_{n-2}\rangle = \langle I_{1}\cap J \rangle \leq \langle I \cap J \rangle\]
by the intersection property of $\Ga_{n}^1$. On the other hand, if $\alpha_{n-1}\not\in I_1$ then 
\[ \lambda_{1}\in \langle I_1 \rangle \cap \langle\alpha_{i+1},\ldots,\alpha_{n-1}\rangle 
\leq \Ga_{n}^2 \cap \Ga_{n}^1 = \langle\alpha_{0},\ldots,\alpha_{n-2}\rangle ,\]
and therefore also 
\[ \lambda_{1}\in \langle\alpha_{0},\ldots,\alpha_{n-2}\rangle \cap \langle\alpha_{i+1},\ldots,\alpha_{n-1}\rangle
= \langle\alpha_{i+1},\ldots,\alpha_{n-2}\rangle = \langle I_{1}\cap J \rangle \leq \langle I \cap J \rangle ,\]
once again by the intersection property of $\Ga_{n}^1$. 

Finally, if $J=\{\alpha_0,\ldots,\alpha_{n-2},\beta_{n-1}\}$ then
\[ \lambda_{1}\in \langle I_1 \rangle \cap \langle J \rangle 
\leq \Ga_{i}^{+} \cap \Ga_{n}^2 = \langle\alpha_{i+1},\ldots,\alpha_{n-2},\beta_{n-1}\rangle. \]
Two possibilities can occur for $I_1$. If $\beta_{n-1}\in I_1$, then $\{\alpha_{i+1},\ldots,\alpha_{n-2},\beta_{n-1}\}= I_1\cap J$ and hence $\lambda_{1}\in\langle I \cap J \rangle$. However, if $\beta_{n-1}\not\in I_1$ then 
\[\lambda_{1}\in \langle I_1 \rangle \cap \langle\alpha_{i+1},\ldots,\alpha_{n-2},\beta_{n-1}\rangle 
\leq \Ga_{n}^1 \cap \Ga_{n}^2 = \langle\alpha_{0},\ldots,\alpha_{n-2}\rangle ,\]
and therefore also
\[\lambda_{1}\in \langle\alpha_{i+1},\ldots,\alpha_{n-2},\beta_{n-1}\rangle \cap \langle\alpha_{0},\ldots,\alpha_{n-2}\rangle
= \langle\alpha_{i+1},\ldots,\alpha_{n-2}\rangle = \langle I_1 \cap J \rangle,\]
now by the intersection property in $\Ga_n^2$. 

Thus in either case $\lambda_1\in\langle I \cap J \rangle$, and since  
$\lambda_0\in\langle I \cap J \rangle$, we also have $\mu\in\langle I \cap J \rangle$. This completes the proof.
\hfill$\square$
\bigskip

\noindent\textbf{Remarks}. 
The intersection conditions explicitly mentioned in Lemma~\ref{intersec} can often be further reduced. For example, when $n=3$,  the three conditions for the mutual intersections of $3$-generator subgroups imply the remaining conditions, provided these subgroups are C-groups.
\medskip

\section{A Universal Construction} Suppose that $\Po$ and $\Q$ are 
regular $n$-polytopes with 
isomorphic facets $\K$ and imagine that we have unlimited copies of each.
Pick a base copy of $\Po$, say, and along each of its facets attach 
a distinct copy of 
$\Q$. Now only certain facets of the new $\Q$'s are exposed, so to them
we attach distinct $\Po$'s. Continue  this in `alternating' fashion.
Presumably we may then
adjoin an infinite number of  copies of $\Po$'s and $\Q$'s
about each ridge of each copy.
But does this process really make sense and do we get an $(n+1)$-polytope 
$\mathcal{S}$ from it? If so, what is $\Gamma(\mathcal{S})$?
Does it matter if instead we begin with a base copy of $\Q$?

We will address these questions by  amalgamating the automorphism
groups of $\Po$ and $\Q$. However, we cannot yet answer a harder question:
does the construction work when around each ridge we attempt to 
alternate just $k$ copies each of $\Po$ and $\Q$ before definitely closing
up, with $k < \infty$?

We begin  by extracting from  \cite[Ch. I, \S 7.3]{bour} some
properties of groups amalgamated 
along subgroups,  as they apply to our construction.
Suppose therefore that we have specified base
flags for $\Po$ and $\Q$ and hence also corresponding lists 
of $n$ standard involutory generators for 
$\Ga(\Po)$ and $\Ga(\Q)$. Then $\Po$ and $\Q$
will have isomorphic facets  precisely when there is an isomorphism
from the facet subgroup of $\Ga(\Po)$ to that of $\Ga(\Q)$
which pairs in order the first $n-1$ generators of $\Ga(\Po)$ 
with those of $\Ga(\Q)$. With these data understood we can unambiguously
amalgamate $\Gamma(\Po)$ with $\Gamma(\Q)$ along
$\Ga(\K)$, giving the group 
$$\Pi : = \Pi(\Po,\Q):= \Ga(\Po) {\ast}_{\Ga(\K)} \Ga(\Q) 
\; .$$ 
But an important feature of the construction is 
that $\Ga(\Po)$ and $\Ga(\Q)$ embed into $\Pi$ 
\cite[Ch.~I, \S 7.3, Prop. 4]{bour}.
Therefore we are justified at the outset in simply assuming that
\begin{eqnarray*}
\Ga(\Po)& =& \langle \alpha_0, \ldots, \alpha_{n-2},\alpha_{n-1}\rangle\,, \\
\Ga(\Q) &= &\langle \alpha_0, \ldots, \alpha_{n-2},\beta_{n-1} \rangle\,, \\
\Ga(\K) &= &\langle \alpha_0, \ldots, \alpha_{n-2}\rangle 
\end{eqnarray*}
are all subgroups of 
$\Pi= \langle \alpha_0, \ldots, \alpha_{n-2},\alpha_{n-1},\beta_{n-1} 
\rangle$. 

One can prove that a presentation of $\Pi$ is obtained by amalgamating
defining relations for $\Ga(\Po)$ (on its generators
$\alpha_0, \ldots, \alpha_{n-2},\alpha_{n-1}$) 
with defining relations for  $\Ga(\Q)$ 
(on $\alpha_0, \ldots, \alpha_{n-2},\beta_{n-1}$).
We can therefore represent this arrangement of groups 
schematically in the diagram

\begin{equation}
\label{unigp}
\centering
\begin{picture}(180,35)
\put(100,0){
\multiput(15,0)(45,0){2}{\circle*{5}}
\multiput(100,-27)(0,53.3){2}{\circle*{5}}
\multiput(-130,0)(45,0){2}{\circle*{5}}
\put(-5,0){\line(1,0){20}}
\put(-85,0){\line(1,0){20}}
\put(-130,0){\line(1,0){45}}
\put(15,0){\line(1,0){45}}
\put(60,0){\line(3,2){40}}
\put(60,0){\line(3,-2){40}}
\put(100,-27){\line(0,1){53.3}}
\put(4,9){\scriptsize $\alpha_{n-3}$}
\put(48,9){\scriptsize $\alpha_{n-2}$}
\put(30,-6){\scriptsize $p_{n-2}$}
\put(-134,9){\scriptsize $\alpha_{0}$}
\put(-89,9){\scriptsize $\alpha_{1}$}
\put(-110,-6){\scriptsize $p_1$}
\put(-49,-0.5){$\ldots\ldots$}
\put(68,21){\scriptsize $p_{n-1}$}
\put(68,-22){\scriptsize $q_{n-1}$}
\put(105,-29){\scriptsize $\beta_{n-1}$}
\put(105,25){\scriptsize $\alpha_{n-1}$}
\put(103.5,-3){\scriptsize $\infty$}}
\end{picture}
\end{equation}
\vskip.4in

\noindent
which evidently has the structure displayed earlier in (\ref{2kgp}).

The group $\Pi$ is usually not a Coxeter group; but it will be
 a quotient of some Coxeter group with a diagram like this. We make no 
 assumptions about  the branch labels $p_1, \ldots, p_{n-2},p_{n-1},q_{n-1}$.
Such a label could 
 even be `$2$', indicating no branch at all.  
 However, it is   a property of free products with amalgamation that
 $\alpha_{n-1} \beta_{n-1}$ has infinite period, which 
 explains the label of the right-hand branch. In other words,
$\langle \alpha_{n-1} ,  \beta_{n-1} \rangle$  is an infinite dihedral group.

For more significant calculations we must employ a standard factorization 
in amalgamated products. Let
$T_{\Po}$ and  $T_{\Q}$  be transversals to 
$\Ga(\K)$ 
in $\Ga(\Po)$ and  $\Ga(\Q)$, respectively,             
with $1 \in T_{\Po} \cap T_{\Q}$. Then every $\mu \in \Pi$
has a unique \textit{reduced decomposition} of length $m \geq 0$, say
\begin{equation}\label{redmu}
\mu = \kappa \tau_1  \cdots \tau_m\; .
\end{equation}
Here $\kappa \in \Ga(\K)$; and  if $m \geq 1$, then all $\tau_i \neq 1$ and
$\tau_i, \tau_{i+1}$ belong to different transversals
$T_{\Po}, T_{\Q}$ for $1 \leq i <n$ \cite[Ch. I, \S 7.3, Prop. 5]{bour}. We call 
$\kappa$  and $\tau_1,\ldots,\tau_m$, respectively, the {\em leading element\/} and the 
{\em transversal elements\/} for the reduced decomposition of $\mu$.

In fact, we shall select our transversals in a special way.

\begin{lemma}\label{utrans}
There are transversals $T_{\Po}$ and $T_{\Q}$ such that for
$0\leq j \leq n-1$, $T_{\Po}$ contains a transversal $T_{\Po, j}$ for
$\langle \alpha_j, \ldots, \alpha_{n-2}\rangle$
in $\langle \alpha_j, \ldots, \alpha_{n-2},\alpha_{n-1} \rangle$;
and $T_{\Q}$ contains a transversal $T_{\Q, j}$ for
$\langle \alpha_j, \ldots, \alpha_{n-2}\rangle$
in $\langle \alpha_j, \ldots, \alpha_{n-2},\beta_{n-1} \rangle$. Moreover, 
\[ \{1,\alpha_{n-1}\}=T_{\Po, n-1}\subseteq T_{\Po, n-2} \subseteq \ldots \subseteq
T_{\Po, 1}\subseteq T_{\Po, 0} = T_{\Po} \]
and 
\[ \{1,\beta_{n-1}\}=T_{\Q, n-1}\subseteq T_{\Q, n-2} \subseteq \ldots \subseteq
T_{\Q, 1}\subseteq T_{\Q, 0} = T_{\Q} .\]
\end{lemma}

\noindent\textbf{Proof}. We proceed by induction on $j$, as $j$ runs from 
$n-1$ down to $0$. When $j = n-1$ we simply pick $T_{\Po, n-1} = \{1, \alpha_{n-1}\}$
as transversal for $\langle 1 \rangle$ in $\langle \alpha_{n-1} \rangle$. 
For the inductive step, suppose that
$\tau_1, \tau_2 \in T_{\Po, j+1}$ (in $\langle \alpha_{j+1}, \ldots, 
\alpha_{n-2},\alpha_{n-1}\rangle$)
satisfy 
$$\tau_1 \equiv \tau_2 \;\mathrm{mod}\;  \langle \alpha_j, \ldots, \alpha_{n-2}\rangle\; .$$
Then 
$$\tau_1^{-1}   \tau_2  \in  \langle \alpha_j, \ldots, \alpha_{n-2}\rangle \cap 
  \langle \alpha_{j+1}, \ldots, \alpha_{n-2},\alpha_{n-1} \rangle\;
  = \langle \alpha_{j+1}, \ldots, \alpha_{n-2}\rangle, $$
by the intersection condition (\ref{interreg}) for string C-groups. Thus $\tau_1 = \tau_2$,
and we may extend $T_{\Po, j+1}$ to $T_{\Po, j}$ for $j = n-2, \ldots, 0$.
\hfill$\square$

\medskip

\begin{lemma}\label{amalj} For $0 \leq j \leq n-1$,
\begin{equation}\label{jstruct}
\langle \alpha_j, \ldots, \alpha_{n-2},\alpha_{n-1},\beta_{n-1} \rangle \simeq
\langle \alpha_j, \ldots, \alpha_{n-2},\alpha_{n-1} \rangle 
\ast_{\langle \alpha_j, \ldots, \alpha_{n-2} \rangle} 
\langle \alpha_j, \ldots, \alpha_{n-2},\beta_{n-1} \rangle\; .
\end{equation}
In particular, 
$\langle \alpha_{n-1},\beta_{n-1} \rangle$ is the infinite dihedral group;
and 
$ \langle \alpha_{n-2},\alpha_{n-1},\beta_{n-1} \rangle$
is the 
\noindent Coxeter group with diagram
\begin{equation}
\centering
\begin{picture}(180,30)
\put(0,0){
\multiput(100,-27)(0,53.3){2}{\circle*{5}}
\put(60,0){\circle*{5}}
\put(60,0){\line(3,2){40}}
\put(60,0){\line(3,-2){40}}
\put(100,-27){\line(0,1){53.3}}
\put(35,-2){\scriptsize $\alpha_{n-2}$}
\put(68,21){\scriptsize $p_{n-1}$}
\put(68,-22){\scriptsize $q_{n-1}$}
\put(105,-29){\scriptsize $\beta_{n-1}$}
\put(105,25){\scriptsize $\alpha_{n-1}$}
\put(103.5,-3){\scriptsize $\infty$}}
\end{picture}
\end{equation}
\vskip.4in
\end{lemma}
\noindent
\textbf{Proof}. For convenience, let 
$A = \langle \alpha_j, \ldots, \alpha_{n-2},\alpha_{n-1} \rangle$,
$B = \langle \alpha_j, \ldots, \alpha_{n-2},\beta_{n-1} \rangle$,\\
$C = \langle \alpha_j, \ldots, \alpha_{n-2} \rangle$, 
all subgroups of $\Pi$.  
We must show that the left side of (\ref{jstruct}) is isomorphic to $A \ast_C B$.
Since the inclusions $A \hookrightarrow \Pi$, 
$B \hookrightarrow \Pi$ agree on $C$, there is a natural map 
$\varphi : A \ast_C B \rightarrow \Pi$ whose image is
$ \langle \alpha_j, \ldots, \alpha_{n-2},\alpha_{n-1},\beta_{n-1} \rangle$.
(We abuse notation a little here.)

Suppose that $\mu \in \ker(\varphi)$
has the reduced decomposition
$\mu = \kappa \tau_1 \cdots \tau_m$, where now 
$\kappa \in C$ and the $\tau_j$'s belong alternately to transversals for $C$ in $A$ or $B$.
By Lemma~\ref{utrans}, these transversals transfer under $\varphi$
to subsets of the transversals $T_{\Po}$, $T_{\Q}$
in $\Pi$. Applying $\varphi$ we therefore get $1 = \kappa \tau_1 \cdots \tau_m$ in 
$\Pi$, where uniqueness of the reduced decomposition  gives $\kappa = 1$ and $m = 0$. 
Thus $\varphi$ is injective. (Compare \cite[Ch. I, Exercise \S 7 - 28]{bour}.)
The   special cases when $j = n-2, n-1$ follow at
once from the standard
presentation of an amalgamated product.
\hfill $\square$

The two previous  lemmas have immediate consequences for the structure of transversals. For $-1\leq j\leq n-2$, let $\Po_j$, $\Q_j$, and $\K_j$ respectively, denote the co-faces of $\Po$, $\Q$, and $\K$, at their basic $j$-face. Then by Lemma~\ref{amalj} the subgroup 
\[ \Pi_{j}^+:=\langle \alpha_{j+1}, \ldots, \alpha_{n-2},\alpha_{n-1},\beta_{n-1} \rangle\] 
of $\Pi$ is isomorphic to the amalgamated product 
\[ \begin{array}{lll}
\Pi(\Po_j,\Q_j) &= & \Ga(\Po_j) \ast_{\Ga(\K_j)} \Ga(\Q_j) \\[.08in]
&=&\langle \alpha_{j+1}, \ldots, \alpha_{n-2},\alpha_{n-1} \rangle 
\ast_{\langle \alpha_{j+1}, \ldots, \alpha_{n-2} \rangle} 
\langle \alpha_{j+1}, \ldots, \alpha_{n-2},\beta_{n-1} \rangle .
\end{array} \]
It follows from the inductive setup of our transversals that, apart from a shift in subscripts, the transversals $T_{\Po_j}$ and $T_{\Q_j}$ associated with $\Po_j$ and $\Q_j$ relative to the amalgamated product $\Pi(\Po_j,\Q_j)$ as in Lemma~\ref{utrans} are just the subsets $T_{\Po,j+1}$ of $T_\Po$ and $T_{\Q,j+1}$ of $T_\Q$. We now have an inductive argument in place to establish the following lemma.

\begin{lemma}\label{indtrans}
Let $T_\Po$ and $T_\Q$ be transversals  inductively built as in 
Lemma~\ref{utrans}, 
let $-1\leq j\leq n-2$, and let 
$\mu=\kappa\tau_{1}\ldots\tau_{m}$ be the reduced decomposition of an 
element $\mu$ in $\Pi_{j}^+$. Then 
$\kappa\in\langle\alpha_{j+1},\ldots,\alpha_{n-2}\rangle$ and 
$\tau_i\in T_{\Po,j+1}\cup T_{\Q,j+1}$ for $i=1,\ldots,m$.
\end{lemma}

\noindent\textbf{Proof}. 
The subgroup $\Pi_{j}^+$ of $\Pi$ is an amalgamated product isomorphic to $\Pi(\Po_j,\Q_j)$, with unique reduced decompositions of its elements $\mu$ relative to the respective transversals $T_{\Po,j+1}=T_{\Po_j}$ and $T_{\Q,j+1}=T_{\Q_{j}}$. Since the latter are subsets of $T_\Po$ and $T_\Q$, respectively, the lemma follows from the uniqueness of the representations of the reduced decompositions in $\Pi(\Po_j,\Q_j)$ and in $\Pi$. 
\hfill$\square$
\medskip

The previous lemmas enable us to establish the crucial intersection property for $\Pi$, required for the construction of the desired semiregular polytope.

\begin{lemma}\label{univinter}
$\Pi$ satisfies the intersection condition \emph{(\ref{intercon})}.
\end{lemma}

\noindent\textbf{Proof}. 
We proceed inductively. By construction, $\Ga(\Po)$ and $\Ga(\Q)$ already are C-groups embedded in $\Pi$. In view of Lemma~\ref{amalj} we also may assume that $\Pi_0$ is a C-group. 

According to Lemma~\ref{intersec} we must verify three intersection conditions, namely $\Ga(\Po)\cap\Ga(\Q)=\Ga(\mathcal{K})$ and both $\Pi_{i}^{+}\cap\Ga(\Po) = \langle \alpha_{i+1}, \ldots,\alpha_{n-1}\rangle$ and $\Pi_{i}^{+}\cap\Ga(\Q)= \langle \alpha_{i+1}, \ldots, \alpha_{n-2},\beta_{n-1}\rangle$ for $i=0,\ldots,n-2$. To this end, suppose $T_\Po$ and $T_\Q$ are transversals chosen as in Lemma~\ref{utrans}. 

First consider the subgroup $\Ga(\Po)\cap\Ga(\Q)\,(=\Pi_n^1\cap \Pi_n^2)$ of $\Pi$. Each element in $\Ga(\Po)$ (resp. $\Ga(\Q)$) has a reduced decomposition with at most one (non-trivial) transversal element from $T_\Po$ and $T_\Q$, respectively. By the uniqueness of reduced decompositions in $\Pi$, an element in $\Ga(\Po)\cap\Ga(\Q)$ cannot involve any transversal elements at all and must lie $\Ga(\K)$. Thus $\Ga(\Po)\cap\Ga(\Q)=\Ga(\K)$. 

Now suppose $\mu\in\Pi_{i}^{+}\cap\Ga(\Po)$ for some $i=0,\ldots,n-2$. By Lemma~\ref{indtrans}, if $\mu=\kappa\tau_{1}\ldots\tau_{m}$ is the reduced decomposition in $\Pi_{i}^+$ (and hence in $\Pi$), then $\kappa\in\langle\alpha_{i+1},\ldots,\alpha_{n-2}\rangle$ and $\tau_j\in T_{\Po,i+1}\cup T_{\Q,i+1}$ for $j=1,\ldots,m$. On the other hand, $\mu\in\Ga(\Po)$, so necessarily $m\leq 1$, and also $\tau_1\in T_{\Po,i+1}$ if $m=1$. In either case $\mu\in \langle \alpha_{i+1}, \ldots,\alpha_{n-2}, \alpha_{n-1}\rangle$.

If $\mu\in\Pi_{i}^{+}\cap\Ga(\Q)$ for some $i=0,\ldots,n-2$, we can similarly conclude that $\mu\in \langle \alpha_{i+1}, \ldots,\alpha_{n-2}, \beta_{n-1}\rangle$, as before with $m\leq 1$ but now with $\tau_1\in T_{\Q,i+1}$ if $m=1$.
\hfill$\square$

\medskip
Now we can apply Theorem~\ref{wytsemi} and Proposition~\ref{2orbit} to the group $\Pi$.

\begin{theorem}\label{univ}
Suppose $\Po$ and $\Q$ are regular $n$-polytopes with isomorphic facets
$\K$. Then the group
$\Pi  = \Ga(\Po) {\ast}_{\Ga(\K)} \Ga(\Q)$ is a group of automorphisms
for a semiregular  $(n+1)$-polytope $\mathcal{U}_{\Po,\Q}$
whose facets are copies of $\Po$ and $\Q$ appearing alternately 
around each  ridge, each a copy of $\K$.
Each section $\mathcal{U}_{\Po,\Q}/\R$ determined  by an $(n-2)$-face $\R$
is an apeirogon $\{\infty\}$. The polytope  $\mathcal{U}_{\Po,\Q}$ is regular
if and only if $\Po$ and $\Q$ are isomorphic. In this case $\Gamma(\mathcal{U}_{\Po,\Q}) =\Pi \rtimes C_2$; otherwise
$\Gamma(\mathcal{U}_{\Po,\Q}) =\Pi$. 
\end{theorem}

\noindent\textbf{Proof}.
It remains to determine the structure of the full automorphism group of 
$\mathcal{U}_{\Po,\Q}$. If $\Po$ and $\Q$ are isomorphic polytopes, 
then $\Pi$ admits an involutory group automorphism that pairs up respective 
generators of $\Ga(\Po)$ and $\Ga(\Q)$ and hence corresponds to the symmetry 
of the underlying diagram for $\Pi$ in its horizontal axis. Note here that 
$\Pi$ admits a presentation that is symmetric with respect to $\Po$ and $\Q$. 
If $\Po$ and $\Q$ are not isomorphic, then clearly such an automorphism cannot 
exist. Now Proposition~\ref{2orbit} applies and completes the proof.

\hfill$\square$ 

In the semiregular polytope $\mathcal{U}_{\Po,\Q}$, an infinite number of 
copies of the given $n$-polytopes $\Po$ and $\Q$ appear alternately as 
facets around each $(n-2)$-face $\R$, determining a section 
$\mathcal{U}_{\Po,\Q}/\R$ of rank $2$ isomorphic to an apeirogon $\{\infty\}$. 
Ideally,  for each integer $k \geq 2$, we would also like to construct 
a similar such semiregular polytope, now with $k$ copies each 
of $\Po$ and $\Q$ appearing alternately around each $(n-2)$-face. At 
the group level, this would involve analysis of the quotient 
$\Pi^{k}=\Pi^{k}(\Po,\Q)$ of
the universal group $\Pi=\Pi(\Po,\Q)$ obtained by adding the extra relation
\begin{equation}
\label{extra}(\alpha_{n-1}\beta_{n-1})^{k}= 1 
\end{equation}
to the defining relations of $\Pi$, that is, by factoring out the normal subgroup $\N^k$ generated by the conjugates of
$(\alpha_{n-1}\beta_{n-1})^{k}$ in $\Pi$. Then three tasks would have to be accomplished, with the first being the most challenging:\ verification of the intersection condition for $\Pi_k$; proof that the groups $\Ga(\Po)$ and $\Ga(\Q)$ really are embedded in $\Pi^k$; and proof that $\alpha_{n-1}\beta_{n-1}$ really has order $k$ in the quotient $\Pi^k$. Solutions for these tasks seem to require rather complicated arguments concerning reduced decompositions of certain elements of the amalgamated product $\Pi$.
\medskip

\section{Polytopes from reflection groups over finite fields}
\label{finfi}

In this section  we briefly sketch a construction of semiregular polytopes 
based on modular reduction techniques applied to certain reflection groups. 
 More details are described in a forthcoming paper~\cite{moscsem}.
 Reflection groups over finite fields and their related regular polytopes 
have been studied in~\cite{monsch1,monsch2,monsch3}.

Suppose that $\Ga = \langle \alpha_0, \ldots, \alpha_{n-2},\alpha_{n-1},\beta_{n-1}\rangle$ 
is an abstract Coxeter group represented, up to branch labels, by a tail-triangle Coxeter
 diagram as in (\ref{2kgp}). Let $p_{i,j}$ and $q_{i,n-1}$ denote the orders of 
$\alpha_{i}\alpha_{j}$ and $\alpha_{i}\beta_{n-1}$ in $\Ga$, so in particular
$p_{i,i}=1$, $p_{j,i} = p_{i,j} =: p_{i+1}$ if $j=i+1$, and $p_{i,j}=2$ otherwise, and 
similarly $q_{n-1,n-1} = k$, $q_{n-1,n-2} = q_{n-2,n-1} =: q_{n-1}$, and $q_{i,n-1}=2$ 
otherwise. (We will not require notation for the order of $\beta_{n-1}^2$.)  Let $V$ be 
real $(n+1)$-space, with basis $\{a_{0},\ldots,a_{n-1},b_{n-1}\}$ and symmetric bilinear 
form $x \cdot y$ defined by
\begin{equation}
\label{bilfor}
a_i \cdot a_j := - 2 \cos{\textstyle \frac{\pi}{p_{i,j}}},\quad\,
a_i \cdot b_{n-1} := -2 \cos{\textstyle\frac{\pi}{q_{i,n-1}}}, \quad\,
b_{n-1} \cdot b_{n-1} := 2.
\end{equation}
Let $R: \Gamma \rightarrow G$ be the standard (faithful) representation of $\Gamma$ in $V$, 
where 
\[ G = \langle r_{0},\ldots,r_{n-1},s_{n-1} \rangle \] 
is the isometric reflection group generated by the reflections with {\em roots\/} 
$a_{0},\ldots,a_{n-1},b_{n-1}$ (see \cite[\S 5.3--5.4]{humph}); thus, 
\[ \begin{array}{rll}
r_i(x) \!\!&=&\!\!x - (x \cdot a_i)\, a_i \quad\; (i = 0,\ldots,n-1), \\[.05in]
s_{n-1}(x)\!\! &=&\!\! x - (x \cdot b_{n-1})\, b_{n-1} .
\end{array} \] 
Then, with respect to the basis $\{a_{0},\ldots,a_{n-1},b_{n-1}\}$ of $V$, 
the reflections $r_{0},\ldots,r_{n-1},s_{n-1}$ are represented by matrices in the general 
linear group $GL_n(\mathbb{D})$ 
over the ring of integers $\mathbb{D}$ in an algebraic number field determined by $\Ga$.
More explicitly, if $\xi$ is a primitive $2m$th 
root of unity, where $m$ denotes the 
lowest common multiple of all finite $p_{i,j}$ and $q_{i,n-1}$, then 
$\mathbb{D} = \mathbb{Z} [\xi]$, the ring of integers in $\mathbb{Q} (\xi)$. Hence we may 
view $G$ as a subgroup of $GL_n(\mathbb{D})$ and {\em reduce $G$ mod $p$}, for any prime 
$p$ (see \cite{monsch1} and~\cite[ch. XII]{cure})). 

The modular reduction technique is most easily described for Coxeter groups $G$ (or $\Ga$) 
which are {\em crystallographic\/}, meaning that $G$ leaves some lattice in $V$ invariant. 
There is a simple combinatorial characterization of crystallographic Coxeter groups, 
which in the present context takes the following form (see \cite{simp}): $G$ is 
crystallographic if and only if $p_{i,j}, q_{i,n-1} = 2$, $3$, $4$, $6$ or 
$\infty$ for all $i,j$ and the triangular circuit of (\ref{2kgp}) (assuming $k\geq 3$) 
contains no, or two, branches labelled $4$ and no, or two, branches labelled $6$. 
If $G$ is crystallographic, then there is a \textit{basic system} 
$\{c_{0},\ldots,c_{n-1},d_{n-1}\}$, with $c_i:=t_i a_i$ and $d_{n-1}=t_{n-1}'b_{n-1}$ 
for certain $t_i > 0$ and $t_{n-1}'>0$, such that 
\[ l_{i,j} := - t_{i}^{-1} (a_{i} \cdot a_{j})\, t_{j} \in \mathbb{Z} \quad (0 \leq i,j 
\leq n-1), \]
\[ m_{i,n-1} := - t_{i}^{-1} (a_{i} \cdot b_{n-1})\, t_{n-1}' \in \mathbb{Z}  
\quad (0 \leq i \leq n-1), \]
and 
\[ m_{n-1, j} := - t_{j} (a_{j} \cdot b_{n-1})\, {t_{n-1}'}^{\!\!\!\!-1} \in \mathbb{Z}  
\quad (0 \leq j \leq n-1). \]
Here, $l_{i,i} =-2$ for all $i$, and $l_{i,j}=0$ if $p_{i,j}=2$; similarly, 
$m_{n-1,n-1} =-2$, and $m_{i,n-1}=0$ if $q_{i,n-1}=2$. 
Then, for the rescaled roots, the generating reflections of $G$ are represented by 
integral matrices and are given by
\begin{equation}
\label{reflnII}
\begin{array}{rll}
r_i(c_j) \!\!&=&\!\! c_j + l_{i,j} c_i , \\[.05in]
r_{i}(d_{n-1})\!\! &=&\!\! d_{n-1} + m_{i,n-1} c_{i},\\[.1in]
s_{n-1}(c_j)\!\! &=&\!\! c_{j} + m_{n-1,j}d_{n-1},\\[.05in]
s_{n-1}(d_{n-1})\!\! &=&\!\! -d_{n-1}.
\end{array}
\end{equation}
The corresponding {\em root lattice} ${(\displaystyle \oplus_{j} \mathbb{Z} c_j) 
\oplus \mathbb{Z} d_{n-1}}$ is $G$-invariant. Thus we may take $\mathbb{D}=\mathbb{Z}$ 
and reduce $G$  modulo 
any positive integer. The most interesting case occurs when 
the reduction is modulo a prime $p$, allowing $p=2$. A crystallographic Coxeter group 
$G$ generally has many basic systems, and the modular reduction may depend 
on the particular  system chosen.

We illustrate our method for the crystallographic Coxeter group
$G$ with diagram 
\begin{equation}
\label{tri1}
\begin{picture}(180,25)
\put(25,0){
\multiput(15,0)(45,0){2}{\circle*{5}}
\multiput(95,-23)(0,46){2}{\circle*{5}}
\put(15,0){\line(1,0){45}}
\put(60,0){\line(3,2){35}}
\put(60,0){\line(3,-2){35}}
\put(74,15){\scriptsize{$4$}}
\put(72,-20){\scriptsize{$\infty$}}
\put(12,7){\scriptsize $r_0$}
\put(56,7){\scriptsize $r_1$}
\put(100,-25.5){\scriptsize $s_2$}
\put(100,21.5){\scriptsize $r_2$}}
\end{picture}
\end{equation}
\vskip.3in 
\noindent
Now $n=3$ and $\{c_0,c_1,c_2,d_2\}$ is a basic system for real $4$-space $V$. As the 
diagram in (\ref{tri1}) has no branches labelled $6$, a change in the underlying basic 
system has little effect on the reduction modulo an odd prime $p$; in fact, passing to 
a different basic system merely results in conjugation inside $GL_4(\mathbb{Z}_p)$. Here 
we may assume that the base vectors $c_0,c_1,c_2,d_2$ have squared length $1,1,2,4$, 
respectively (see \cite[\S 4]{monsch1}). 

When $p$ is odd, the reduced group $G^p$ is generated by four isometric reflections 
in the finite orthogonal geometry determined on $V^p:=\mathbb{Z}_{p}^4$ by the 
reduction of the corresponding (integral) bilinear form modulo $p$. However, $V^p$ 
is non-singular only when $p>3$, since its discriminant is $-6 \bmod{p}$. 
Thus the cases  $p=2$ and $p=3$ require special treatment.

Now, by Lemma~\ref{intersec} and the subsequent remarks, $G^p$ is a tail-triangle 
C-group if and only if the $3$-generator subgroups 
$\langle r_0,r_1,r_2 \rangle^p$, $\langle r_0,r_1,s_2 \rangle^p$ and 
$\langle r_1,r_2,s_2 \rangle^p$ are C-groups satisfying  
\begin{equation}
\label{threeint}
\begin{array}{rcl} 
\langle r_0,r_1,r_2 \rangle^p \cap \langle r_0,r_1,s_2\rangle^p &=& \langle r_0,r_1\rangle^p,\\
\langle r_1,r_2,s_2 \rangle^p \cap \langle r_0,r_1,r_2\rangle^p &=& \langle r_1,r_2\rangle^p,\\
\langle r_1,r_2,s_2 \rangle^p \cap \langle r_0,r_1,s_2\rangle^p &=& \langle r_1,s_2\rangle^p.\\
\end{array} 
\end{equation}

When $p$ is odd, each of the three $3$-dimensional subspaces of $V^p$ 
spanned by $\{c_0,c_1,c_2\}$, $\{c_0,c_1,d_2\}$ or $\{c_1,c_2,d_2\}$, respectively, 
is left invariant by its respective $3$-generator subgroup, and is non-singular, 
with discriminant $\tfrac{1}{2}$, $-1$ or $-4 \bmod{p}$. Moreover, each of the 
three $2$-dimensional subspaces of $V^p$ spanned by $\{c_0,c_1\}$, $\{c_1,c_2\}$ 
or $\{c_1,d_2\}$, respectively, is invariant under the respective dihedral subgroup 
$\langle r_0,r_1\rangle^p$, $\langle r_1,r_2\rangle^p$ or $\langle r_1,s_2\rangle^p$, 
and has discriminant $\tfrac{3}{4}$, $1$ or $0\bmod{p}$; here the third subspace is 
singular (for all $p$), but the first and second subspaces are non-singular, except 
for the first when $p=3$. Now, without elaborating in detail, it turns out that the 
above $3$-generator subgroups of $G^p$ really are C-groups, and that the three 
intersection conditions in (\ref{threeint}) can be verified using methods similar to 
those described in \cite{monsch1,monsch2} (when $p\geq 5$); more precisely, the first 
two intersections conditions require an analog of \cite[Theorem 4.2]{monsch1}, and the 
third an
analog of \cite[Corollary 3.2]{monsch2}. The case $p=3$ can be verified by hand, or 
using GAP~\cite{gap}.   

Now Theorem~\ref{wytsemi} and Proposition~\ref{2orbit} tell us that $G^p$ is the 
automorphism group of a semiregular $4$-polytope $\mathcal{S}$ whose two 
kinds of facets $\Po$ and $\Q$ are regular polyhedra determined by the 
subgroups $\langle r_0,r_1,r_2 \rangle^p$ and $\langle r_0,r_1,s_2 \rangle^p$. 
In particular, $\Po$ is the octahedron $\{3,4\}$; and $\Q$ is the regular map 
of type $\{3,p\}$ of \cite[\S 5.7]{monsch1}, with automorphism group 
$G^p = O_{1}(3,p,0)$ if $p\geq 5$, and $G^p\simeq \mathbb{S}_4$ if $p=3$. 
(Recall that $O_{1}(3,p,0)$ denotes the subgroup generated by the reflections 
of spinor norm $1$ in the full orthogonal group $O(3,p,0)$ of a non-singular 
orthogonal space $\mathbb{Z}_p^3$. Note that 
$O_{1}(3,p,0) \simeq PSL_{2}(\mathbb{Z}_p)\rtimes C_2$ and $O(3,p,0) \simeq 
PGL_{2}(\mathbb{Z}_p)\rtimes C_2$.)  
When $p=3$, $5$ or $7$, respectively, $\Q$ is the tetrahedron $\{3,3\}$, 
icosahedron $\{3,5\}$ or Klein map $\{3,7\}_8$. Each edge of $\mathcal{S}$ 
lies in four facets, namely two copies of $\Po$ and two copies of $\Q$, 
occurring alternately. The automorphism group $G^p$ of $\mathcal{S}$ is 
given by 
\begin{equation}\label{autes}
\Gamma({\cal S}) = 
\left\{ \begin{array}{ll}
O_1(4,p,1)\;, & \mbox{ if } p \equiv 1,7 \bmod{24},\\
O(4,p,1)\;, & \mbox{ if } p \equiv 5,11 \bmod{24},\\
O_1(4,p,-1)\;, & \mbox{ if } p \equiv 17,23 \bmod{24},\\
O(4,p,-1)\;, & \mbox{ if } p \equiv 13,19 \bmod{24} ,
\end{array}\right. \; .
\end{equation}
(Recall that $O_{1}(4,p,\varepsilon)$, with $\epsilon=\pm 1$, denotes the subgroup 
generated by the reflections of spinor norm $1$, in the full orthogonal group 
$O(4,p,\epsilon)$ of a non-singular orthogonal space $\mathbb{Z}_p^4$. Here 
$\varepsilon=1$ if the Witt index is $2$, and $\varepsilon=-1$ if the Witt index is $1$.)

For $p=3$ the space $V^3$ is singular with $1$-dimensional radical. The group $G^3$ 
has order $1296$ and consists of all isometries which fix the radical pointwise, so that
$G^3 \simeq \mathbb{Z}_3^3 \rtimes \langle r_0, r_1, r_2 \rangle^3$. Thus we can think of
$G^3$ as a crystallographic group with finite invariant lattice
$\mathbb{Z}_3^3$ and octahedral point group 
$\langle r_0, r_1, r_2 \rangle^3 \simeq O(3,3,0) \simeq B_3$. Referring to
\cite{bouw}, it is now easy to see why $G^3$ is isomorphic to the full automorphism group
of the Gray graph (see also \cite{semi}).

We now consider the special prime $p=2$.  Of the several essentially distinct 
basic systems admitted by $G$, only that in which the rescaled base vectors 
have squared lengths $1,1,2,4$ actually yields a tail-triangle C-group.
We find that $G^2$ has order $96$ and, a little unexpectedly, that 
$r_1 s_2$ has period $4$. The semiregular polytope $\mathcal{S} = \mathcal{S}(G^2)$
has two sets of $4$ hemioctahedral facets distributed among a meager $3$ vertices. 
Each vertex-figure is a toroid $\{4,4\}_{(2,0)}$. Since $G^2$ does admit the automorphism
mentioned in Proposition~\ref{2orbit}(a), $\mathcal{S}$ is in fact regular; and  since the facets
are combinatorially flat, $\mathcal{S}$ must coincide with the (flat) universal regular polytope 
$\{\{3,4\}_3,\{4,4\}_{(2,0)}\}$ (see \cite[4E5]{arp}) denoted  $\{3,4,4\}\ast 192a$ in the Census~\cite{harcen}.  
It is clear from our comment in Example~\ref{tomeg} that $G^2$ is isomorphic to the 
automorphism group 
$\Gamma(\mathcal{T}_{hh}) = \langle \rho_0, \rho_1, \rho_2, \rho_3\rangle$ for the tomotope.
This isomorphism is induced by mapping
$(r_0, r_1, r_2, s_2)$ to $(\rho_0, \rho_1, \rho_3, \rho_0\rho_2)$, so the geometrical 
connection between $\mathcal{S}$ and $\mathcal{T}_{hh}$ is a little obscure.

The diagram in (\ref{tri1}) can be viewed in three ways as a tail-triangle diagram 
for the same underlying group $G$ (or $G^p$). The corresponding intersection conditions 
(\ref{threeint}) involve the same $3$-generator subgroups in each case. Thus $G^p$ is a 
tail-triangle C-group in three different ways, and hence can give 
rise to three mutually 
non-isomorphic semiregular $4$-polytopes. We have already discussed the semiregular 
$4$-polytope $\mathcal{S}$ associated with the original diagram.

Now when the diagram is taken in the form 
\begin{equation}
\label{tri2}
\begin{picture}(180,25)
\put(25,0){
\multiput(15,0)(45,0){2}{\circle*{5}}
\multiput(95,-23)(0,46){2}{\circle*{5}}
\put(15,0){\line(1,0){45}}
\put(60,0){\line(3,2){35}}
\put(60,0){\line(3,-2){35}}
\put(35,-8){\scriptsize{$4$}}
\put(72,-20){\scriptsize{$\infty$}}}
\end{picture}
\end{equation}
\vskip.3in
\noindent
the corresponding semiregular $4$-polytope $\mathcal{S}'$ again has two kinds of facets 
$\Po$ and $\Q$, occurring alternately around every edge of $\mathcal{S}'$. Now 
$\Po$ is the $3$-cube $\{4,3\}$; and $\Q$ is the regular map of type $\{4,p\}$ of 
\cite[\S 5.9]{monsch1}, with group $G^p = O_{1}(3,p,0)$ if $p\equiv \pm 1\bmod{8}$, and 
$G^p=O(3,p,0)$ otherwise. When $p\geq 5$ the automorphism group of $\mathcal{S}'$ is 
again the group in (\ref{autes}). However, when $p=3$ the $4$-polytope $\mathcal{S}'$ 
is regular and is isomorphic to the 
regular $3$-toroid $\{4,3,4\}_{(3,3,0)}$ (see \cite[6D]{arp}). 

Finally, from  
\begin{equation}  
\label{tri3}
\begin{picture}(180,25)
\put(25,0){
\multiput(15,00)(45,0){2}{\circle*{5}}
\multiput(95,-23)(0,46){2}{\circle*{5}}
\put(15,0){\line(1,0){45}}
\put(60,0){\line(3,2){35}}
\put(60,0){\line(3,-2){35}}
\put(35,-8){\scriptsize{$\infty$}}
\put(72,15){\scriptsize{$4$}}}
\end{picture}
\end{equation} 
\vskip.3in
\noindent
we obtain a semiregular $4$-polytope $\mathcal{S}''$ whose  facets $\Po$ of type 
$\{p,4\}$ are the duals of those of type $\{4,p\}$ of $\mathcal{S}'$, and
whose facets $\Q$ of type $\{p,3\}$ are the duals of those of type $\{3,p\}$ of $\mathcal{S}$. 
When $p=3$, we can view $\mathcal{S}''$ as a semiregular tessellation of the $3$-torus by 
$54$ tetrahedra  and $27$ octahedra. 
This $3$-torus can be obtained by identifying opposite faces of the parallelepiped
spanned by vectors $(3,3,0), (3,0,3)$ and $(0,3,3)$ (see Figure~\ref{fig1}). 
When $p\geq 5$, the automorphism  group of $\mathcal{S}''$ is again the group in (\ref{autes}).

The finite semiregular polytopes $\mathcal{S}$, $\mathcal{S}'$ and $\mathcal{S}''$ derived 
from the diagrams for $G^p$ in (\ref{tri1}), (\ref{tri2}) and (\ref{tri3}) are quotients of 
the infinite semiregular polytopes associated with the infinite Coxeter groups with these diagrams as Coxeter diagrams. In fact, the obvious group epimorphisms mapping generators to generators determine
coverings between the infinite polytopes and the corresponding finite polytopes 
(see \cite[2D]{arp}). 
\vskip.15in

\noindent\textbf{Acknowledgements}.
We are grateful to Gordon Williams for help with the figures. We also want to thank
the referee for various useful comments.

\end{document}